\documentclass[sn-mathphys,Numbered]{sn-jnl}%

\usepackage{graphicx}%
\usepackage{multirow}%
\usepackage{amsmath,amssymb,amsfonts}%
\usepackage{amsthm}%
\usepackage{mathrsfs}%
\usepackage[title]{appendix}%
\usepackage{xcolor}%
\usepackage{textcomp}%
\usepackage{manyfoot}%
\usepackage{booktabs}%
\usepackage{algorithm}%
\usepackage{algorithmicx}%
\usepackage{algpseudocode}%
\usepackage{enumitem}
\usepackage{listings}%
\usepackage{tikz}
\usetikzlibrary{shapes,calc,arrows,positioning,backgrounds,patterns}
\usepackage{comment}
\theoremstyle{thmstyleone}%
\newtheorem{theorem}{Theorem}%
\newtheorem{proposition}[theorem]{Proposition}%
\newtheorem{assumption}[theorem]{Assumption}

\theoremstyle{thmstyletwo}%
\newtheorem{example}{Example}%
\newtheorem{remark}{Remark}%
\newtheorem{lemma}{Lemma}
\newtheorem{notation}{Notation}
\newtheorem{corollary}{Corollary}
\theoremstyle{thmstylethree}%
\newtheorem{definition}{Definition}%

\newcommand{\W}[1]{{#1}}

\makeatletter %
\def\@acknow{}%
\long\def\EarlyAcknow#1 \par{%
\def\@acknow{\abstractfont\abstracthead*{Acknowledgments}%
#1\par}}%

\def\printabstract{\ifx\@acknow\empty\else\@acknow\fi\par%
    \ifx\@abstract\empty\else\@abstract\fi\par}
\makeatother
\begin{document}

\title[The pressure-wired Stokes element:
  	a mesh-robust version of the Scott-Vogelius element]{The pressure-wired Stokes element:
  	a mesh-robust version of the Scott-Vogelius element}

\author[1]{\fnm{Benedikt} \sur{Gr\"{a}\ss le}}\email{graesslb@math.hu-berlin.de}

\author[2]{\fnm{Nis-Erik} \sur{Bohne}}\email{nis-erik.bohne@uzh.ch}

\author*[2]{\fnm{Stefan} \sur{Sauter}}\email{stas@math.uzh.ch}

\affil[1]{\orgdiv{Institut f\"{u}r
		Mathematik}, \orgname{Humboldt-Universit\"{a}t zu Berlin}, \orgaddress{\street{Rudower Chaussee 25}, \city{Berlin}, \postcode{12489}, \state{Berlin}, \country{Germany}}}

\affil[2]{\orgdiv{Institut f\"{u}r
		Mathematik}, \orgname{Universit\"{a}t Z\"{u}rich}, \orgaddress{\street{Winterthurerstrasse 190}, \city{Z\"{u}rich}, \postcode{CH-8057}, \state{Z\"{u}rich}, \country{Switzerland}}}

		\EarlyAcknow{
The first author is supported by the \emph{Deutsche Forschungsgemeinschaft} (DFG) under Germany's Excellence Strategy -- The Berlin Mathematics Research Center MATH+ (EXC-2046/1, project ID: 390685689).}
\abstract{The Scott-Vogelius finite element pair for the numerical discretization of
the stationary Stokes equation in 2D is a popular element which is based on
a continuous velocity approximation of polynomial order $k$ and a
discontinuous pressure approximation of order $k-1$. It employs a \textquotedblleft singular
distance\textquotedblright\ (measured by some geometric mesh quantity $
\Theta \left( \mathbf{z}\right) \geq 0$ for triangle vertices $\mathbf{z}$)
and imposes a local side condition on the pressure space associated to
vertices $\mathbf{z}$ with $\Theta \left( \mathbf{z}\right) =0$. The method
is inf-sup stable for any fixed regular triangulation and $k\geq 4$.
However, the inf-sup constant deteriorates if the triangulation contains
nearly singular vertices $0<\Theta \left( \mathbf{z}\right) \ll 1$. 

In this paper, we introduce a very simple parameter-dependent modification
of the Scott-Vogelius element such that the inf-sup constant is independent
of nearly-singular vertices. We will show by analysis and also by
numerical experiments that the effect on the divergence-free condition for
the discrete velocity is negligibly small.}

\keywords{$hp$ finite elements, Scott-Vogelius
	elements, inf-sup stability, mass conservation}

\pacs[MSC Classification]{65N30, 65N12, 76D07}

\maketitle

\section{Introduction}\label{sec1}

In this paper we consider the numerical solution of the
stationary Stokes equation in a bounded two-dimensional polygonal Lipschitz
domain by a conforming Galerkin finite element method. 

\textbf{Motivation.} The intuitive choice for a Stokes element $\left( 
\mathbf{S}_{k}\left( \mathcal{T}\right) ,\mathbb{P}_{k-1}\left( \mathcal{T}%
\right) \right) $, where $\mathbf{S}_{k}\left( \mathcal{T}\right) $ denotes
the space of continuous velocity fields with local polynomial degree $k$
over a given triangulation $\mathcal{T}$ and $\mathbb{P}_{k-1}\left( 
\mathcal{T}\right) $ the discontinuous pressures space of degree $k-1$, is
in general not inf-sup stable (see, e.g., \cite{vogelius1983right} and{\Huge %
\ }\cite[Chap. 7]{Braessengl} for quadrilateral meshes). A careful analysis 
\cite{vogelius1983right}, \cite{ScottVogelius} of the range of the
divergence operator reveals that the \W{dimension of the} image $\operatorname{div}(\mathbf{S}_{k}(%
\mathcal{T}))\subset \mathbb{P}_{k-1}(\mathcal{T})$ reduces by one linear
constraint for every \emph{singular vertex}. For $k\geq 4$, this results in
the inf-sup stable Scott-Vogelius \cite{ScottVogelius} pair $(\mathbf{S}_{k}(%
\mathcal{T}),\operatorname{div}(\mathbf{S}_{k}(\mathcal{T})))$ on families of
shape-regular meshes \cite{GuzmanScott2019} that naturally computes fully
divergence-free velocity approximations. However the Scott-Vogelius element, as described in \cite{vogelius1983right}%
, \cite{ScottVogelius} and \cite{GuzmanScott2019}, has two major drawbacks.

\begin{enumerate}
\item The inf-sup constant is not robust with respect to small perturbations
of the mesh when a singular vertex becomes nearly singular, measured through the geometric quantity $\Theta_{\min} :=
\W{\min}\left\{\left. \Theta \left(\mathbf{z}\right) \; \right\vert \; \Theta \left(\mathbf{z} \right) > 0\right\}$. Here $%
\Theta \left( \mathbf{z}\right) $ is a measure of the \textit{singular
distance} of a vertex $\mathbf{z}$ from a proper singular situation. $\Theta_{\min}$ strongly
affects the stability of the discretization, the size of the
discretization error, as well as the condition number of the resulting
algebraic linear system. 
\item In order to implement the method, the condition: \textquotedblleft Is
a mesh point, say $\mathbf{z}$\textbf{,} a singular
vertex?\textquotedblright\ cannot be realized in floating point arithmetics
and has to be replaced by a threshold condition for \textquotedblleft nearly
singular vertices\textquotedblright\ of the form \textquotedblleft $\Theta
\left( \mathbf{z}\right) \leq \varepsilon $\textquotedblright. %
Through this threshold condition, it is possible that the constraints are imposed on nearly singular vertices, i.e. for
$ \W{0}<\Theta \left( \mathbf{z}\right) \leq \varepsilon $ and the discrete velocity looses the divergence-free property as a consequence. To the best of our knowledge, this effect has not been analyzed in the literature and we will
estimate the influence of $\varepsilon >0$ to the divergence-free property
of the discrete velocity in Section \ref{SecSolFree}. 
\end{enumerate}
In \cite[Rem.~2]{AP_Locking}, two mesh modification strategies are sketched
as a remedy of (1): i) nearly singular vertices are moved so that
they become properly singular; ii) triangles which contain a nearly singular
vertex are refined. Both strategies require the finite element code to have
control on the mesh 
\W{generator. 
	However, in many engineering finite element codes the mesh design is
separated from the discretization and hence, interactive refinement features are not
available.
In addition fast solvers in computational fluid dynamics sometimes make use of a very regular
mesh structure (e.g., C-meshes for the modelling of airfoils) and then local
mesh refinement is not compatible.}

In this paper, we propose a simpler strategy where a modification of the
mesh is avoided. We introduce a parameter dependent modification to the
standard Scott-Vogelius pair, the \emph{pressure-wired Stokes
element}, and prove that the inf-sup constant is independent of $\Theta
_{\min }$, the mesh width, the polynomial degree but depends only on the
shape-regularity of the mesh. 

\textbf{Main Contributions.} The construction of our modification employs a
geometric quantity $\Theta \left( \mathbf{z}\right) \geq 0$ which measures
the \textit{singular distance} of the vertex $\mathbf{z}$ in the mesh from a
singular configuration. For some arbitrary (in general small) control
parameter $0\leq \eta $, the pressure space is reduced by the same
constraint as in the Scott-Vogelius FEM for every vertex with singular
distance $\Theta (\mathbf{z})\leq \eta $, which we call \emph{nearly singular%
}. Since the constraints involve the \W{alternating} sum of pressure values over all
triangles that encircle a nearly singular vertex \W{(similar to a wire that is loosely woven around the tips of all triangles
sharing a common vertex)} we call the new element the 
\textit{pressure-wired Stokes element}. The proof that the inf-sup constant of this element is independent of the
mesh width, the polynomial degree, and depends on the mesh only \W{through} the
shape-regularity constant will be based on the results in \cite[Sections 4
\& 5]{AP_Locking}, where the discrete stability is proved via a lower bound for the inf-sup constant of the form $%
c\Theta _{\min } >0$ with a positive constant $c$ only
depending on the domain and the shape regularity. %
We can therefore choose the threshold $\eta$ for our
generalization of the Scott-Vogelius element such
that the resulting pressure\W{-}wired Stokes element is
inf-sup stable independent of the mesh size, the
polynomial degree, and any (nearly) singular vertex. These improvements
come at some cost: in general we cannot expect divergence-free velocity
approximations as for the standard Scott-Vogelius pair. 
We investigate the dependence of the $L^{2}$ norm of the divergence of the
velocity approximation $\mathbf{u}_{\mathbf{S}}$ in the pressure-wired
Stokes element and establish an estimate of the form $\left\Vert \operatorname{div}%
\mathbf{u}_{\mathbf{S}}\right\Vert _{L^{2}\left( \Omega \right) }\leq C\eta $
with $C>0$ being independent of the mesh width and the polynomial degree. We
emphasize that this estimate does not require any regularity of the Stokes
\W{solution}. Numerical experiments will be reported in Section \ref{sec:Numerical experiments} and show that the constant $C$ is very small for the considered examples.

\textbf{Literature overview.} 
\W{The Scott-Vogelius pair for triangulations of two-dimensional domains (see \cite{vogelius1983right},
	\cite{ScottVogelius}) is inf-sup stable for $k\geq4$, while
	the inf-sup constant deteriorates for nearly singular vertices.
	This problem can be attenuated by using special mesh-refinement strategies, e.g., Alfeld splits that are also known as
	barycentric refinements and provide inf-sup stability already for $k\geq 2$ in two dimensions
	\cite{AQ:QuadraticVelocityLinear1992} and
	$k\geq 3$ in three dimensions \cite{Zha:NewFamilyStable2005}.
In three dimensions, inf-sup stability and
a full characterisation of the divergence
$\operatorname*{div}\mathbf{S}_{k}(\mathcal{T})$ is not known on general meshes, see \cite{FMS:TwoConjecturesStokes2022} and references
therein.
}
Our pressure-wired Stokes element
places few additional constraints on the pressure space and acquires robust
inf-sup stability on arbitrary shape-regular grids -- the proof for the
estimate of the inf-sup constant is based on the theory developed in \cite%
{vogelius1983right}, \cite{ScottVogelius}, \cite%
{GuzmanScott2019}, \cite{AP_Locking}, \cite{Bernardi_Raugel_1985}, \cite{SauterCR_prob}. 
An alternative approach constitutes the enrichment of the velocity space,
e.g., in \cite{merdon_SV} with Raviart-Thomas bubble functions leading to a
divergence-free velocity approximation in $H(\operatorname{div})$. Conforming
alternatives to the Scott-Vogelius element are the Mini element \cite%
{arnold1984stable}, the (modified) Taylor-Hood element \cite[Chap. 3, \S 7]%
{Braessengl}, the Bernardi-Raugel element \cite{Bernardi_Raugel_1985}, and
the element by Falk-Neilan \cite{FalkNeilan} to mention some but few of
them.
We note that there exist further possibilities to obtain a stable
discretization of the Stokes equation; one is the use of non-conforming
schemes and/or modifications of the discrete equation by adding stabilizing
terms. We do not go into details here but refer to the monographs and
overviews \cite{BoffiBrezziFortin}, \cite{ErnGuermondII}, \cite%
{Brenner_Crouzeix} instead.
\W{For a discussion on the divergence-free condition and pressure-robustness, we refer to \cite{JLM+:DivergenceConstraintMixed2017}.
}

\textbf{Further contributions and outline:} After introducing the
Stokes problem on the continuous as well as on the discrete level and the
relevant notation in Section \ref{NuDiscrete}, we define the conforming
pressure-wired Stokes element in Section \ref{SecRSV}. The first main result
in Section \ref{SecProof} establishes the discrete inf-sup condition for
this new element with a lower bound on the inf-sup constant that is
independent of $h,k$, and (nearly) critical points. The second main result
controls the $L^{2}$ norm of the divergence of the discrete velocity in
Section \ref{SecSolFree} and verifies its negligibility for small $\eta \ll
1 $ in practice. In fact, the $L^{2}$ norm tends to zero (at least) linearly
in $\eta $ without imposing any regularity assumption on the continuous
problem. The involved constants are again independent of $h$, $k$, and $%
\Theta _{\min }$ but possibly depend on the mesh via its shape-regularity
constant. We report on numerical evidence on optimal convergence rates in\
Section \ref{sec:Numerical experiments}, both as an $h$-version with regular
mesh refinement and as a $k$-version that successively increases the
polynomial degree. While this new element is technically not
divergence-free, our benchmarks suggest that the discrete divergence is 
\textit{near machine-precision} already on coarse meshes and for moderate
parameters $\eta >0$. %

\section{The Stokes problem and its numerical discretization\label%
{NuDiscrete}}

Let $\Omega \subset \mathbb{R}^{2}$ denote a bounded polygonal Lipschitz
domain with boundary $\partial \Omega $. We consider the numerical solution
of the Stokes equation%
\[
\begin{array}{lll}
-\Delta \mathbf{u}+\nabla p & =\mathbf{f} & \text{in }\Omega , \\ 
\operatorname{div}\mathbf{u} & =0 & \text{in }\Omega%
\end{array}%
\]%
with homogeneous Dirichlet boundary conditions for the velocity and the
usual normalisation condition for the pressure%
\[
\mathbf{u}=\mathbf{0}\quad \text{on }\partial \Omega \quad \text{and\quad }%
\int_{\Omega }p=0. 
\]%
Throughout this paper standard notation for Lebesgue and Sobolev spaces
applies. All function spaces are considered over the field of real numbers. The space $%
H_{0}^{1}\left( \Omega \right) $ is the closure of the space of smooth functions with compact support in $\Omega $ with respect
to the $H^{1}\left( \Omega \right) $ norm. Its dual space is denoted by $%
H^{-1}\left( \Omega \right) :=H_{0}^{1}\left( \Omega \right) ^{\prime }$.
The scalar product and norm in $L^{2}\left( \Omega \right) $ read%
\[
\left( u,v\right) _{L^{2}\left( \Omega \right) }:=\int_{\Omega }uv\quad 
\text{and}\quad \left\Vert u\right\Vert _{L^{2}\left( \Omega \right)
}:=\left( u,u\right) _{L^{2}\left( \Omega \right) }^{1/2}\qquad \text{in }%
L^{2}\left( \Omega \right) .%
\]%
Vector-valued and $2\times 2$ tensor-valued analogues of the function spaces
are denoted by bold and blackboard bold letters, e.g., $\mathbf{H}^{s}\left(
\Omega \right) =\left( H^{s}\left( \Omega \right) \right) ^{2}$ and $\mathbb{%
H}^{s}\left( \Omega \right) =\left( H^{s}\left( \Omega \right) \right)
^{2\times 2}$ and analogously for other quantities.

\begin{notation}
	\label{NotVectors}
	For vectors $\mathbf{v},\mathbf{w}\in\mathbb{R}^{2}$, the
	Euclidean scalar product $\left\langle \mathbf{v},\mathbf{w}%
	\right\rangle :=\sum_{j=1}^{2}v_{j}w_{j}$ induces the Euclidean norm 
	$\left\Vert \mathbf{v}\right\Vert :=\left\langle \mathbf{v},\mathbf{v}%
	\right\rangle ^{1/2}$. 
	We denote the area of a subset $M\subset \mathbb R^2$ by 
	$\left\vert M\right\vert $  and
	write \W{$[  \mathbf{v}|\mathbf{w}]  $} for
	the $2\times2$ matrix with column vectors $\mathbf{v},\mathbf{w}\W{\in\mathbb R^2}$.
\end{notation}
The $\mathbf{L}^{2}\left( \Omega\right) $ scalar product and norm for vector
valued functions are%
\[
\left( \mathbf{u},\mathbf{v}\right) _{\mathbf{L}^{2}\left( \Omega\right)
}:=\int_{\Omega}\left\langle \mathbf{u},\mathbf{v}\right\rangle \quad \text{%
and\quad}\left\Vert \mathbf{u}\right\Vert _{\mathbf{L}^{2}\left(
\Omega\right) }:=\left( \mathbf{u},\mathbf{u}\right) _{\mathbf{L} ^{2}\left(
\Omega\right) }^{1/2}. 
\]
In a similar fashion, we define for $\mathbf{G},\mathbf{H}\in\mathbb{L}%
^{2}\left( \Omega\right) $ the scalar product and norm by%
\[
\left( \mathbf{G},\mathbf{H}\right) _{\mathbb{L}^{2}\left( \Omega\right)
}:=\int_{\Omega}\left\langle \mathbf{G},\mathbf{H}\right\rangle \quad \text{%
and\quad}\left\Vert \mathbf{G}\right\Vert _{\mathbb{L}^{2}\left(
\Omega\right) }:=\left( \mathbf{G},\mathbf{G}\right) _{\mathbb{L} ^{2}\left(
\Omega\right) }^{1/2}, 
\]
where $\left\langle \mathbf{G},\mathbf{H}\right\rangle
=\sum_{i,j=1}^{2}G_{i,j}H_{i,j}$. Finally, let $L_{0}^{2}\left(
\Omega\right) :=\left\{ u\in L^{2}\left( \Omega\right)
:\int_{\Omega}u=0\right\} $. We introduce the bilinear forms $a:\mathbf{H}%
^{1}\left( \Omega\right) \times\mathbf{H}^{1}\left( \Omega\right) \rightarrow%
\mathbb{R}$ and $b: \mathbf{H}^{1}\left( \Omega\right) \times L^2
\left(\Omega\right) \to \mathbb{R}$ by%
\begin{equation}
\begin{array}{ccc}
a\left( \mathbf{u},\mathbf{v}\right) :=\left( \nabla\mathbf{u} ,\nabla%
\mathbf{v}\right) _{\mathbb{L}^{2}\left( \Omega\right) } & \text{ and } & b
\left(\mathbf{u},p\right) = \left(\operatorname{div} \mathbf{u},
p\right)_{L^2\left(\Omega\right)}, \label{defabili}%
\end{array}%
\end{equation}
where $\nabla\mathbf{u}$ and $\nabla\mathbf{v}$ denote the gradients of $%
\mathbf{u}$ and $\mathbf{v}$. Given $\mathbf{F}\in\mathbf{H}%
^{-1}\left( \Omega\right)$, the variational form of the stationary Stokes
problem seeks $\left( \mathbf{u},p\right) \in\mathbf{H}_{0}^{1}\left(
\Omega\right) \times L_{0}^{2}\left( \Omega\right)$ such that 
\begin{equation}
\begin{array}{lll}
a\left( \mathbf{u},\mathbf{v}\right) -b\left( \mathbf{v},p\right) & =\mathbf{%
F}\left( \mathbf{v}\right) & \forall\mathbf{v}\in\mathbf{H}_{0}^{1}\left(
\Omega\right) , \\ 
b \left( \mathbf{u},q\right) & =0 & \forall q\in L_{0}^{2}\left(
\Omega\right) .%
\end{array}
\label{varproblemstokes}
\end{equation}
The inf-sup condition guarantees well-posedness of \eqref{varproblemstokes},
cf. \cite{Girault86} for details. In this paper, we consider a conforming
Galerkin discretization of (\ref{varproblemstokes}) by a pair $\left( 
\mathbf{S},M\right) $ of finite dimensional subspaces of the continuous
solution spaces $\left( \mathbf{H}_{0}^{1}\left( \Omega\right)
,L_{0}^{2}\left( \Omega\right) \right) $. For any given $\mathbf{F}\in%
\mathbf{H}^{-1}\left( \Omega\right) $ the weak formulation yields $\left( 
\mathbf{u}_{\mathbf{S}},p_{M}\right) \in\mathbf{S}\times M\;$such that%
\begin{equation}
\begin{array}{lll}
a\left( \mathbf{u}_{\mathbf{S}},\mathbf{v}\right) -b\left( \mathbf{v}
,p_{M}\right) & =\mathbf{F}\left( \mathbf{v}\right) & \forall\mathbf{v} \in%
\mathbf{S}, \\ 
b\left( \mathbf{u}_{\mathbf{S}},q\right) & =0 & \forall q\in M.%
\end{array}
\label{discrStokes}
\end{equation}
It is well known that the bilinear form $a\left( \cdot,\cdot\right) $ is
symmetric, continuous, and coercive so that problem (\ref{discrStokes}) is
well-posed if the bilinear form $b\left( \cdot,\cdot\right) $ satisfies the
inf-sup condition.

\begin{definition}
	Let $\mathbf{S}$ and $M$ be finite-dimensional subspaces of $\mathbf{H}%
	_{0}^{1}\left(  \Omega\right)  $ and $L_{0}^{2}\left(  \Omega\right)  $. The
	pair $\left(  \mathbf{S},M\right)  $ is \emph{inf-sup stable} if the \emph{inf-sup constant} $\beta \left(
	\mathbf{S},M\right)  $ is positive, i.e.,%
	\begin{equation} \beta \left(  \mathbf{S},M\right):=
		\inf_{q\in M\backslash\left\{  0\right\}  }\sup_{\mathbf{v}\in\mathbf{S}%
			\backslash\left\{  \mathbf{0}\right\}  }\frac{\left(  q,\operatorname*{div}%
			\mathbf{v}\right)  _{L^{2}\left(  \Omega\right)  }}{\left\Vert \mathbf{v}%
			\right\Vert _{\mathbf{H}^{1}\left(  \Omega\right)  }\left\Vert q\right\Vert
			_{L^{2}\left(  \Omega\right)  }}  >0.
		\label{infsupcond}%
	\end{equation}
	
\end{definition}

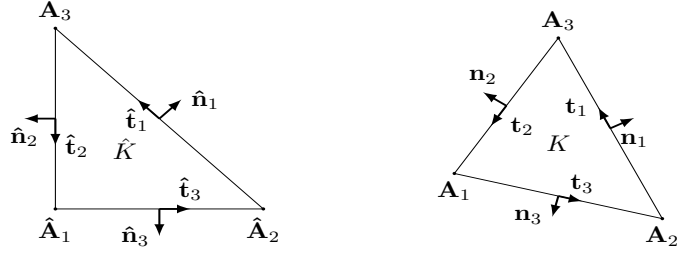
\begin{figure}[]
\centering
\resizebox{.7\textwidth}{10em}{
	\begin{tikzpicture}[scale=.7]
		\node (mc 1) at (0,0) {};
		\node (mc 2) at (4,0) {};
		\node (mc 3) at (0,4) {};
		\node (ms 1) at (2,2) {};
		\node (ms 2) at (0,2) {};
		\node (ms 3) at (2,0) {};
		\node (T) at (1.333, 1.333) {$\hat K$};
		\draw (mc 1.center) -- (mc 2.center) -- (mc 3.center) -- cycle;
		\fill [black] (mc 1) circle (1pt) node[below] {$\mathbf{\hat{A}}_{1}$};
		\fill [black] (mc 2) circle (1pt) node[below] {$\mathbf{\hat{A}}_{2}$};
		\fill [black] (mc 3) circle (1pt) node[above] {$\mathbf{\hat{A}}_{3}$};
		\draw [-latex, thick] (ms 1.center) -- ($(ms 1)!0.2121!90:(mc 2)$) node[right] {$\mathbf{\hat n}_1$};
		\draw [-latex, thick] (ms 2.center) -- ($(ms 2)!0.3!90:(mc 3)$) node[below] {$\mathbf{\hat n}_2$};
		\draw [-latex, thick] (ms 3.center) -- ($(ms 3)!0.3!90:(mc 1)$) node[left] {$\mathbf{\hat n}_3$};
		\draw [-latex, thick] (ms 1.center) -- ($(ms 1)!0.2121!180:(mc 2)$) node[below] {$\mathbf{\hat t}_1$};
		\draw [-latex, thick] (ms 2.center) -- ($(ms 2)!0.3!180:(mc 3)$) node[right] {$\mathbf{\hat t}_2$};
		\draw [-latex, thick] (ms 3.center) -- ($(ms 3)!0.3!180:(mc 1)$) node[above] {$\mathbf{\hat t}_3$};
		\phantom{\draw[black] (mc 1) circle (.2cm);}
	\end{tikzpicture}\hspace{5em}
	\begin{tikzpicture}[scale=.7]
		\node (mc 1) at (0,1) {};
		\node (mc 2) at (4,0) {};
		\node (mc 3) at (2,4) {};
		\node (ms 1) at (3,2) {};
		\node (ms 2) at (1,2.5) {};
		\node (ms 3) at (2,0.5) {};
		\node (T) at (2, 1.666) {$K$};
		\draw (mc 1.center) -- (mc 2.center) -- (mc 3.center) -- cycle;
		\fill [black] (mc 1) circle (1pt) node[below] {$\mathbf{{A}}_{1}$};
		\fill [black] (mc 2) circle (1pt) node[below] {$\mathbf{{A}}_{2}$};
		\fill [black] (mc 3) circle (1pt) node[above] {$\mathbf{{A}}_{3}$};
		\draw [-latex, thick] (ms 1.center) -- ($(ms 1)!0.23!90:(mc 2)$) node[yshift=-1em] {$\mathbf{n}_1$};
		\draw [-latex, thick] (ms 2.center) -- ($(ms 2)!0.3!90:(mc 3)$) node[above] {$\mathbf{ n}_2$};
		\draw [-latex, thick] (ms 3.center) -- ($(ms 3)!0.23!90:(mc 1)$) node[left] {$\mathbf{ n}_3$};
		\draw [-latex, thick] (ms 1.center) -- ($(ms 1)!0.23!180:(mc 2)$) node[left] {$\mathbf{t}_1$};
		\draw [-latex, thick] (ms 2.center) -- ($(ms 2)!0.3!180:(mc 3)$) node[xshift=1.2em] {$\mathbf{t}_2$};
		\draw [-latex, thick] (ms 3.center) -- ($(ms 3)!0.23!180:(mc 1)$) node[above] {$\mathbf{t}_3$};
		\phantom{\draw[black] (mc 1) circle (.2cm);}
\end{tikzpicture}}
\caption{Reference triangle (left) and physical triangle (right)}
\label{fig:T_notation}
\end{figure}

\section{The pressure-wired Stokes element\label{SecRSV}}

Throughout this paper, $\mathcal{T}%
$ denotes a conforming triangulation of the bounded polygonal Lipschitz domain $\Omega
\subset \mathbb{R}^{2}$ into closed triangles: the intersection of two
different triangles is either empty, a common edge, or a common point. The
set of edges in $\mathcal{T}$ is denoted by $\mathcal{E}\left( \mathcal{T}%
\right) $, comprised of boundary edges $\mathcal{E}_{\partial\Omega}\left(  \mathcal{T}\right)  :=\left\{  E\in\mathcal{E}\left( \mathcal{T}\right)  :E\subset\partial \Omega \right\}  $ and interior edges $\mathcal{E}_{\Omega}\left(  \mathcal{T}\right)  :=\mathcal{E}\left( \mathcal{T}\right)  \backslash\mathcal{E}_{\partial\Omega}\left( \mathcal{T}\right)  $. 
For any edge $E\in \mathcal{E}(\mathcal{T})$, we fix a unit normal vector $%
\mathbf{n}_{E}$ with the convention that for boundary edges
 $E\in\mathcal{E}_{\partial\Omega}, \ 
\mathbf{n}_{E}$ points to the exterior of $\Omega $. The set of vertices in 
$\mathcal{T}$ is denoted by $\mathcal{V}\left( \mathcal{T}\right) $ while
the subset of boundary vertices is $\mathcal{V}_{\partial \Omega }\left( 
\mathcal{T}\right) :=\left\{ \mathbf{z}\in \mathcal{V}\left( \mathcal{T}%
\right) :\mathbf{z}\in \partial \Omega \right\} $. The interior vertices
form the set $\mathcal{V}_{\Omega }\left( \mathcal{T}\right) :=\mathcal{V}%
\left( \mathcal{T}\right) \backslash \mathcal{V}_{\partial \Omega }\left( 
\mathcal{T}\right) $. For \W{a triangle} $K\in \mathcal{T}$ \W{with diameter $h_K:= \operatorname*{diam}K$}, the set of its vertices is
denoted by $\mathcal{V}\left( K\right) $. For $\mathbf{z}\in \mathcal{V}%
\left( \mathcal{T}\right) $, we consider the local element vertex patch%
\begin{equation}
\begin{array}{ll}
\mathcal{T}_{%
\mathbf{z}}:=\left\{ K\in \mathcal{T}\mid \mathbf{z}\in K\right\}  & \text{%
and\quad }\omega _{\mathbf{z}}:=\bigcup\limits_{K\in \mathcal{T}_{\mathbf{z}%
}}K%
\end{array}
\label{nodalpatch}
\end{equation}%
with the local mesh width 
$h_{\mathbf{z}}:=\max \left\{ h_{K}:K\in \mathcal{T}_{\mathbf{z}}\right\} $.
For any vertex $\mathbf{z}\in \mathcal{V}\left( \mathcal{T}\right) $, we fix
a counterclockwise numbering of the $N_{\mathbf{z}}:=\left\vert \mathcal{T}_{%
\mathbf{z}}\right\vert $ triangles in 
\begin{equation}
\mathcal{T}_{\mathbf{z}}=\left\{ \left. K_{j} \; \right\vert \; 1\leq j\leq N_{\mathbf{z}}\right\}
\quad \text{and set\ }\mathcal{E}_{\mathbf{z}}:=\left\{ E\in \mathcal{E}%
\left( \mathcal{T}\right) :\mathbf{z}\text{ is an endpoint of }E\right\} .
\label{eqn:enum_Tz}
\end{equation}%
\begin{figure}[tbp]
\center
\W{\resizebox{\textwidth}{8em}{
	\begin{tikzpicture}[scale=0.8]
				\draw (0,0) node[] (z) {$\mathbf{z}$};
		\coordinate (A) at (4,0);
		\coordinate (B) at ($(z)!0.5!120:(A)$);
		\coordinate (C) at ($(z)!0.9!180:(A)$);
		\coordinate (D) at ($(z)!0.53!300:(A)$);
		\draw (A) -- (B) -- (C) -- (D) -- cycle;
		\draw (z) -- (A) node[midway,above] {$E_1$};
		\draw (z) -- (B) node[midway,right] {$E_2$};
		\draw (z) -- (C) node[midway,above] {$E_3$}; 
		\draw (z) -- (D) node[midway,left] {$E_4$};
		\node (T1) at (barycentric cs:z=1,A=1,B=1) {$K_1$};
		\node (T1) at (barycentric cs:z=1,C=1,B=2) {$K_2$};
		\node (T1) at (barycentric cs:z=1,A=1,D=1) {$K_4$};
		\node (T1) at (barycentric cs:z=1,C=1,D=1) {$K_3$};
	\end{tikzpicture}
	\begin{tikzpicture}[scale=0.8]
				\draw (0,0) node[] (z) {$\mathbf{z}$};
		\coordinate (A) at (4,0);
		\coordinate (B) at ($(z)!0.5!120:(A)$);
		\coordinate (C) at ($(z)!0.9!180:(A)$);
		\coordinate (D) at ($(z)!0.53!300:(A)$);
		\draw (A) -- (B) -- (C) -- (D);
		\draw (z) -- (A) node[midway,above] {$E_1$};
		\draw (z) -- (B) node[midway,right] {$E_2$};
		\draw (z) -- (C) node[midway,above] {$E_3$}; 
		\draw (z) -- (D) node[midway,left] {$E_4$};
		\node (T1) at (barycentric cs:z=1,A=1,B=1) {$K_1$};
		\node (T1) at (barycentric cs:z=1,C=1,B=2) {$K_2$};
		\node (T1) at (barycentric cs:z=1,C=1,D=1) {$K_3$};
		\fill[pattern=north east lines, ] ([xshift=0.35em]z.south) -- (D) -- ([xshift=.3em]D) -- ([xshift=.65em]z.south);
		\fill[pattern=north east lines, ] (z.east) rectangle ([yshift=-.3em]A);
	\end{tikzpicture}
	\begin{tikzpicture}[scale=0.8]
				\draw (0,0) node[] (z) {$\mathbf{z}$};
		\coordinate (A) at (4,0);
		\coordinate (B) at ($(z)!0.5!120:(A)$);
		\coordinate (C) at ($(z)!0.9!180:(A)$);
		\coordinate (D) at ($(z)!0.53!300:(A)$);
		\draw (A) -- (B) -- (C);		\draw (z) -- (A) node[midway,above] {$E_1$};
		\draw (z) -- (B) node[midway,right] {$E_2$};
		\draw (z) -- (C) node[midway,above] {$E_3$}; 
		\node (T1) at (barycentric cs:z=1,A=1,B=1) {$K_1$};
		\node (T1) at (barycentric cs:z=1,C=1,B=2) {$K_2$};
		\phantom{		\node (T1) at (barycentric cs:z=1,A=1,D=1) {$K_4$};
		\node (T1) at (barycentric cs:z=1,C=1,D=1) {$K_3$};}
		\fill[pattern=north east lines, ] (z.east) rectangle ([yshift=-.3em]A);
		\fill[pattern=north east lines, ] (z.west) rectangle ([yshift=-.3em]C);
	\end{tikzpicture}
	\begin{tikzpicture}[scale=0.8]
				\draw (0,0) node[] (z) {$\mathbf{z}$};
		\coordinate (A) at (4,0);
		\coordinate (B) at ($(z)!0.5!80:(A)$);
		\coordinate (C) at ($(z)!0.9!180:(A)$);
		\coordinate (D) at ($(z)!0.53!300:(A)$);
		\draw (A) -- (B);		\draw (z) -- (A) node[midway,above] {$E_1$};
		\draw (z) -- (B) node[midway,right] {$E_2$};
		\phantom{		\node (T1) at (barycentric cs:z=1,A=1,D=1) {$K_4$};
		\node (T1) at (barycentric cs:z=1,C=1,D=1) {$K_3$};}
		\node (T1) at (barycentric cs:z=1,A=1,B=1) {$K_1$};
		\fill[pattern=north east lines, ] (z.east) rectangle ([yshift=-.3em]A);
		\fill[pattern=north west lines, ] ([xshift=0.08em]z.north) -- (B) -- ([xshift=-.3em]B) -- ([xshift=-.22em]z.north);
\end{tikzpicture}
}}
\resizebox{\textwidth}{8em}{
	\begin{tikzpicture}[scale=0.8]
				\draw (0,0) node[] (z) {$\mathbf{z}$};
		\coordinate (A) at (4,0);
		\coordinate (B) at ($(z)!0.5!80:(A)$);
		\coordinate (C) at ($(z)!0.9!200:(A)$);
		\coordinate (D) at ($(z)!0.53!300:(A)$);
		\draw (A) -- (B) -- (C) -- (D) -- cycle;
		\draw (z) -- (A) node[midway,above] {$E_1$};
		\draw (z) -- (B) node[midway,right] {$E_2$};
		\draw (z) -- (C) node[midway,above] {$E_3$}; 
		\draw (z) -- (D) node[midway,left] {$E_4$};
		\node (T1) at (barycentric cs:z=1,A=1,B=1) {$K_1$};
		\node (T1) at (barycentric cs:z=1,C=1,B=1) {$K_2$};
		\node (T1) at (barycentric cs:z=1,A=1,D=1) {$K_4$};
		\node (T1) at (barycentric cs:z=1,C=1,D=1) {$K_3$};
	\end{tikzpicture}
	\begin{tikzpicture}[scale=0.8]
				\draw (0,0) node[] (z) {$\mathbf{z}$};
		\coordinate (A) at (4,0);
		\coordinate (B) at ($(z)!0.5!80:(A)$);
		\coordinate (C) at ($(z)!0.9!200:(A)$);
		\coordinate (D) at ($(z)!0.53!300:(A)$);
		\draw (A) -- (B) -- (C) -- (D);
		\draw (z) -- (A) node[midway,above] {$E_1$};
		\draw (z) -- (B) node[midway,right] {$E_2$};
		\draw (z) -- (C) node[midway,above] {$E_3$}; 
		\draw (z) -- (D) node[midway,left] {$E_4$};
		\node (T1) at (barycentric cs:z=1,A=1,B=1) {$K_1$};
		\node (T1) at (barycentric cs:z=1,C=1,B=1) {$K_2$};
		\node (T1) at (barycentric cs:z=1,C=1,D=1) {$K_3$};
		\fill[pattern=north east lines, ] ([xshift=0.35em]z.south) -- (D) -- ([xshift=.3em]D) -- ([xshift=.65em]z.south);
		\fill[pattern=north east lines, ] (z.east) rectangle ([yshift=-.3em]A);
	\end{tikzpicture}
	\begin{tikzpicture}[scale=0.8]
				\draw (0,0) node[] (z) {$\mathbf{z}$};
		\coordinate (A) at (4,0);
		\coordinate (B) at ($(z)!0.5!80:(A)$);
		\coordinate (C) at ($(z)!0.9!200:(A)$);
		\coordinate (D) at ($(z)!0.53!300:(A)$);
		\draw (A) -- (B) -- (C);		\draw (z) -- (A) node[midway,above] {$E_1$};
		\draw (z) -- (B) node[midway,right] {$E_2$};
		\draw (z) -- (C) node[midway,above] {$E_3$}; 
		\node (T1) at (barycentric cs:z=1,A=1,B=1) {$K_1$};
		\node (T1) at (barycentric cs:z=1,C=1,B=1) {$K_2$};
		\phantom{		\node (T1) at (barycentric cs:z=1,A=1,D=1) {$K_4$};
		\node (T1) at (barycentric cs:z=1,C=1,D=1) {$K_3$};}
		\fill[pattern=north west lines, ] ([yshift=-0.25em]z.west) -- (C) -- ([yshift=-.3em]C) -- ([yshift=-.55em]z.west);
		\fill[pattern=north east lines, ] (z.east) rectangle ([yshift=-.3em]A);
	\end{tikzpicture}
	\begin{tikzpicture}[scale=0.8]
				\draw (0,0) node[] (z) {$\mathbf{z}$};
		\coordinate (A) at (4,0);
		\coordinate (B) at ($(z)!0.5!80:(A)$);
		\coordinate (C) at ($(z)!0.9!200:(A)$);
		\coordinate (D) at ($(z)!0.53!300:(A)$);
		\draw (A) -- (B);		\draw (z) -- (A) node[midway,above] {$E_1$};
		\draw (z) -- (B) node[midway,right] {$E_2$};
		\phantom{		\node (T1) at (barycentric cs:z=1,A=1,D=1) {$K_4$};
		\node (T1) at (barycentric cs:z=1,C=1,D=1) {$K_3$};}
		\node (T1) at (barycentric cs:z=1,A=1,B=1) {$K_1$};
		\fill[pattern=north east lines, ] (z.east) rectangle ([yshift=-.3em]A);
		\fill[pattern=north west lines, ] ([xshift=0.08em]z.north) -- (B) -- ([xshift=-.3em]B) -- ([xshift=-.22em]z.north);
\end{tikzpicture}}
\caption{Singular configuration (top) and generic configuration (bottom) of a vertex patch for an interior vertex $\mathbf{z}\in \mathcal{V}%
_{\Omega }(\mathcal{T})$ with $N_{\mathbf{z}}=4$ (resp.~boundery vertex $%
\mathbf{z}\in \mathcal{V}_{\partial \Omega }(\mathcal{T})$ with $N_{\mathbf{z}}=1,2,3$)
triangles}
\label{fig:vertex_patch}
\end{figure}
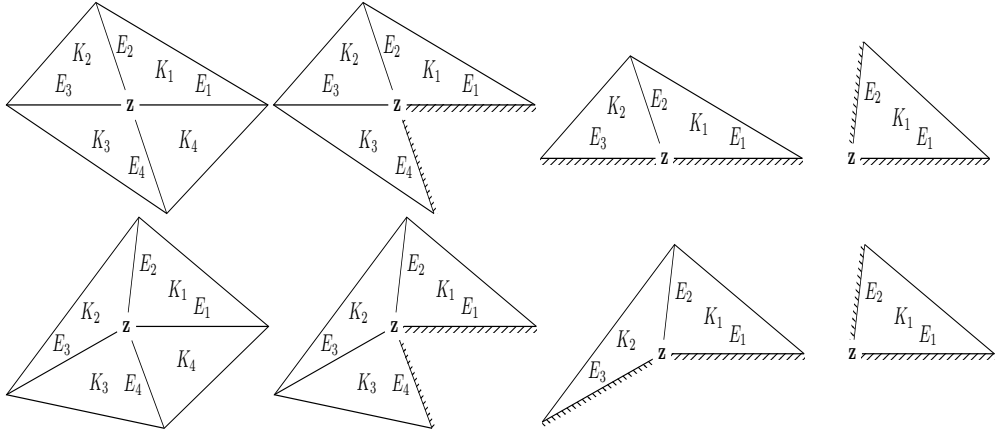
In Fig. \ref{fig:vertex_patch}, $\mathcal{T}_{\mathbf{z}}$ is shown for four important configurations.
The shape-regularity constant 
\begin{equation}
\gamma _{\mathcal{T}}:=\max_{K\in \mathcal{T}}\frac{h_{K}}{\rho _{K}}
\label{defgammat}
\end{equation}%
relates the local mesh width $h_{K}\W{=}\operatorname{diam}K$ with the diameter $%
\rho _{K}$ of the largest inscribed ball in an element $K\in \mathcal{T}$.
The global mesh width is given by $h_{\mathcal{T}}:=\max \left\{ h_{K}:K\in 
\mathcal{T}\right\} $. 
\begin{remark}
	\label{Remangle}The shape-regularity implies the existence of some $\phi_{\mathcal{T}}>0$
	exclusively depending on $\gamma_{\mathcal{T}}$
	such that every triangle angle in $\mathcal{T}$ is bounded from below by
	$\phi_{\mathcal{T}}$. In turn, every triangle angle in $\mathcal{T}$ is
	bounded from above by $\pi-2\phi_{\mathcal{T}}$.
\end{remark}
For $\mathbf{z}\in \mathcal{V}_{\partial \Omega }\left( \mathcal{%
T}\right) $, denote by $\alpha _{\mathbf{z}}$ the interior angle between the
two edges in $\mathcal{E}_{\partial \Omega }\left( \mathcal{T}\right) $ with
joint endpoint $\mathbf{z}$ and regarded from the exterior complement $\mathbb{R}%
^{2}\backslash \overline{\Omega }$. Let%
\begin{equation}
\alpha _{\mathcal{T}}:=\min_{\mathbf{z}\in \mathcal{V}_{\partial \Omega
}\left( \mathcal{T}\right) }\alpha _{\mathbf{z}}  \label{defalphat}
\end{equation}%
denote the minimal outer angle at the boundary vertices that lies between $%
0<\alpha _{\mathcal{T}}<2\pi $ for the Lipschitz domain $\Omega $. 
Let $\mathbb{P}%
_{k}(K)$ denote the space of polynomials up to degree $k\in \mathbb{N}_{0}$
defined on $K\in \mathcal{T}$ and define 
\begin{equation}
\begin{aligned} \mathbb{P}_{k}\left( \mathcal{T}\right) &:=\left\{ q\in
L^{2}\left( \Omega\right) \mid\forall K\in\mathcal{T}:\left. q\right\vert
_{\overset{\circ}{K}}\in\mathbb{P}_{k}\left(\overset{\circ}{K}\right) \right\} ,\\ \mathbb{P}_{k,0}\left(
\mathcal{T}\right) &:=\left\{ q\in\mathbb{P}_{k}\left( \mathcal{T}\right)
\mid\int _{\Omega}q=0\right\}\equiv \mathbb P_k(\mathcal T)\cap
L^2_0(\Omega) ,\\ S_{k}\left( \mathcal{T}\right) &:=\mathbb{P}_{k}\left(
\mathcal{T}\right) \cap H^{1}(\Omega),\\ S_{k,0}\left( \mathcal{T}\right)
&:=S_{k}\left( \mathcal{T} \right) \cap H_{0}^{1}\left(\Omega\right) .
\end{aligned}  \label{Pkdefs}
\end{equation}%
The vector-valued spaces are $\mathbf{S}_{k}\left( \mathcal{T}\right)
:=S_{k}\left( \mathcal{T}\right) ^{2}$ and $\mathbf{S}_{k,0}\left( \mathcal{T%
}\right) :=S_{k,0}\left( \mathcal{T}\right) ^{2}$. It is well known that the
most intuitive Stokes element $\left( \mathbf{S}_{k,0}\left( \mathcal{T}%
\right) ,\mathbb{P}_{k-1,0}\left( \mathcal{T}\right) \right) $ is in general
unstable. %
The analysis in \cite{vogelius1983right}, \cite{ScottVogelius} for $k\geq 4$
relates the instability of $\left( \mathbf{S}_{k,0}\left( \mathcal{T}\right)
,\mathbb{P}_{k-1,0}\left( \mathcal{T}\right) \right) $ to \emph{critical} or 
\emph{singular points} of the mesh $\mathcal{T}$. The set of $\eta 
$-critical points $\mathcal{C}_{\mathcal{T}}\left( \eta \right) $ for the
control parameter $\eta \geq 0$ recovers the definition of the \emph{%
classical critical points} $\mathcal{C}_{\mathcal{T}}(0)$ (introduced in 
\cite[R.1, R.2]{vogelius1983right} and called \emph{\ singular vertices} in 
\cite{vogelius1983right}) for $\eta =0$.

\begin{definition}
\label{DefCritpoint} The local measure of singularity $\Theta \left( \mathbf{%
z}\right) $ at $\mathbf{z}\in \mathcal{V}(\mathcal{T})$ is given by 
\begin{equation}
\Theta \left( \mathbf{z}\right) :=\begin{cases} \max\left\{ \left.
\left\vert \sin\left( \theta_{i}+\theta_{i+1}\right) \right\vert \;
\right\vert \; 0 \leq i \leq N_{\mathbf{z}} \right\} & \text{if
}\mathbf{z}\in\mathcal{V}_{\Omega}\left( \mathcal{T} \right) ,\\ \max\left\{
\left. \left\vert \sin\left( \theta_{i}+\theta_{i+1}\right) \right\vert \;
\right\vert \; 0 \leq i \leq N_{\mathbf{z} -1}\right\} & \text{if
}\mathbf{z}\in\mathcal{V}_{\partial\Omega}(\mathcal T)\wedge
N_{\mathbf{z}}>1,\\ 0 & \text{if
}\mathbf{z\in}\mathcal{V}_{\partial\Omega}(\mathcal T)\wedge
N_{\mathbf{z}}=1, \end{cases}  \label{defthetaz}
\end{equation}%
where the angles $\theta _{j}$ in $\mathcal{T}_{z}$ are numbered
counterclockwise from $1\leq j\leq N_{\mathbf{z}}$ (see (\ref{eqn:enum_Tz}))
and cyclic numbering is applied, i.e. $\theta _{N_{\mathbf{z}+1}}=\theta _{1}
$. For $\eta \geq 0$, the vertex $\mathbf{z}\in \mathcal{V}\left( \mathcal{T}%
\right) $ is an $\eta $\emph{-critical} vertex if $\Theta \left( \mathbf{z}%
\right) \leq \eta $. The set of all $\eta $-critical vertices is given by 
\[
\mathcal{C}_{\mathcal{T}}\left( \eta \right) :=\left\{ \mathbf{z}\in 
\mathcal{V}\left( \mathcal{T}\right) \mid \Theta \left( \mathbf{z}\right)
\leq \eta \right\} .
\]%
For a vertex $\mathbf{z}\in \mathcal{V}\left( \mathcal{T}\right) $, the
functional $A_{\mathcal{T},\mathbf{z}}:\mathbb{P}_{k}\left( \mathcal{T}%
\right) \rightarrow \mathbb{R}$ alternates counterclockwise through the
numbered triangles $K_{\ell }$, $1\leq \ell \leq N_{\mathbf{z}}$ in the patch $%
\mathcal{T}_{\mathbf{z}}$, and is given by%
\begin{equation}
A_{\mathcal{T},\mathbf{z}}\left( q\right) :=\sum_{\ell =1}^{N_{\mathbf{z}}}\left(
-1\right) ^{\ell }\left( \left. q\right\vert _{K_{\ell }}\right) \left( 
\mathbf{z}\right) .  \label{DefATz}
\end{equation}
\end{definition}
\W{From Definition \eqref{defthetaz} it follows that for any $\eta\geq1$ we have
$\mathcal{C}_{\mathcal{T}}\left(  \eta\right)  =\mathcal{V}\left(
\mathcal{T}\right)  $ and the assumption $\eta\in\left[  0,1\right]  $
throughout this paper is not an actual restriction. }
\W{The first row in Figure~\ref{fig:vertex_patch} displays all four singular mesh configurations (i.e.,
$\Theta(\mathbf{z})=0$) that occurs when all edges of the vertex patch $\mathcal{T}(\mathbf{z})$ lie on two straight
lines.
It follows \W{from} \cite{SauterCR_prob} that for small $\eta$, the only possible $\eta$-critical mesh configurations are
perturbations of those four situations and
illustrated in the second row of Figure~\ref{fig:vertex_patch}.
\begin{lemma}[{\cite[Lem.~2.13]{SauterCR_prob}}]\label{lem:eta_config}
There exists a constant $\eta_{0}>0$ only depending on the
shape-regularity of the mesh and the minimal outer angle $\alpha_{\mathcal{T}%
}$ of $\Omega$ such that \W{$\Theta(\mathbf{z})\leq \eta_{0}$} for $\mathbf{z}\in\mathcal{V}(\mathcal{T})$ implies
\begin{enumerate}[label=(\alph*)]
	\item if $\mathbf{z}\in \mathcal{V}_{\Omega}(\mathcal{T})$ is an interior vertex, then $N_{\mathbf{z}}=4$,
	\item if $\mathbf{z}\in \mathcal{V}_{\partial\Omega}(\mathcal{T})$ is a boundary vertex, then
		$N_{\mathbf{z}}\leq 3$.\qed
\end{enumerate}
\end{lemma}
}
The pressure space \W{of the pressure-wired Stokes element is obtained by} requiring the condition $A_{\mathcal{T},%
\mathbf{z}}\left( q\right) =0$ for all $\eta $-critical points.

\begin{definition}
	\label{DefAz}For $\eta\geq0$, the subspace $M_{\eta,k-1}\left(  \mathcal{T}%
	\right) \subset \mathbb{P}_{k-1,0}\left(  \mathcal{T}%
	\right)  $ of the pressure space is given by
	\begin{equation}
		M_{\eta,k-1}\left(  \mathcal{T}\right)  :=\left\{  q\in\mathbb{P}%
		_{k-1,0}\left(  \mathcal{T}\right)  \mid\forall\mathbf{z}\in\mathcal{C}%
		_{\mathcal{T}}\left(  \eta\right)  :A_{\mathcal{T},\mathbf{z}}\left(
		q\right)  =0\right\}  . \label{defMSV}%
	\end{equation}
	The \emph{pressure-wired Stokes element} is given by $\left(  \mathbf{S}%
	_{k,0}\left(  \mathcal{T}\right)  ,M_{\eta,k-1}\left(  \mathcal{T}\right)  \right)  $.
\end{definition}
Note that for the choice $\eta = 0$, $M_{0,k-1}\left( \mathcal{T}\right) $
is the pressure space introduced by Vogelius \cite{vogelius1983right} and
Scott-Vogelius \cite{ScottVogelius} and the following inclusions hold: for $%
0\leq\eta\leq \eta^{\prime}$%
\begin{equation}
M_{\eta^{\prime},k-1}\left( \mathcal{T}\right) \subset M_{\eta,k-1}\left( 
\mathcal{T}\right) \subset M_{0,k-1}\left( \mathcal{T}\right) =Q_{h}^{k-1}
\label{press_incl}
\end{equation}
(with the pressure space $Q_{h}^{k-1}$ in \cite[p. 517]{GuzmanScott2019}).
The existence of a continuous right-inverse of the divergence operator $%
\operatorname{div}: \mathbf{S}_{k,0}(\mathcal{T})\to M_{0,k-1}(\mathcal{T})$ was
proved in \cite{vogelius1983right} and \cite[Thm. 5.1]{ScottVogelius}.

\begin{proposition}
	[Scott-Vogelius]Let $k\geq4$. For any $q\in M_{0,k-1}\left(  \mathcal{T}%
	\right)  $ there exists some $\mathbf{v}\in\mathbf{S}_{k,0}\left(
	\mathcal{T}\right)  $ such that%
	\[
	\operatorname*{div}\mathbf{v}=q\mathbf{\quad}\text{and\quad}\left\Vert
	\mathbf{v}\right\Vert _{\mathbf{H}^{1}\left(  \Omega\right)  }\leq C\left\Vert
	q\right\Vert _{L^{2}\left(  \Omega\right)  }.
	\]
The constant $C$ is independent of $h_{\mathcal{T}}$ and only depends on the shape-regularity of the mesh, the polynomial degree $k$, and on $\Theta_{\min}$, where
	\begin{equation}
		\Theta_{\min}:=\min_{\mathbf{z}\in\mathcal{V}\left(  \mathcal{T}\right)
			\backslash\mathcal{C}_{\mathcal{T}}\left(  0\right)  }\Theta\left(
		\mathbf{z}\right)  . \label{thetamin}%
	\end{equation}
\end{proposition}It follows from \cite[Thm. 1]{GuzmanScott2019} that the
inf-sup constant of the Scott-Vogelius element $\left( \mathbf{S}%
_{k,0}\left( \mathcal{T}\right) ,M_{0,k-1}\left( \mathcal{T}\right) \right) $
can be bounded from below by $\beta \left( \mathbf{S}_{k,0}\left( \mathcal{T}%
\right) ,M_{0,k-1}\left( \mathcal{T}\right) \right) \geq \Theta _{\min
}c\left( \gamma _{\mathcal{T}},k\right) $, where $c\left( \gamma _{\mathcal{T%
}},k\right) >0$ depends on the shape regularity constant $\gamma _{\mathcal{T%
}}$ and the polynomial degree $k$. The dependence on $k$ has been analysed
in a series of papers starting from \cite[Lem. 2.5]{vogelius1983right} with
the final result in \cite[Theorem 5.1]{AP_Locking}, which states that 
\[
\beta \left( \mathbf{S}_{k,0}\left( \mathcal{T}\right) ,M_{0,k-1}\left( 
\mathcal{T}\right) \right) \geq \Theta _{\min }c\left( \gamma _{\mathcal{T}%
}\right) >0
\]%
is independent of the polynomial degree $k$. However, the estimate is \textit{%
non-robust} with respect to small perturbations of critical configurations $%
0\neq \Theta \left( \mathbf{z}\right) \ll 1$. Our notion of $\eta $-critical
points can be regarded as a robust generalization and we will analyse the
consequences in this paper. In the next section, we will prove for the
pressure-wired Stokes element that%
\[
\beta \left( \mathbf{S}_{k,0}\left( \mathcal{T}\right) ,M_{\eta ,k-1}\left( 
\mathcal{T}\right) \right) \geq \left( \Theta _{\min }+\eta \right) c\left(
\gamma _{\mathcal{T}}\right) >0
\]%
by modifying the arguments in \cite[Sections 4 \& 5]{AP_Locking} and
investigate the effect on the divergence-free property of the discrete
velocity in Section \ref{SecSolFree}.

\section{Inf-sup stability of the pressure-wired Stokes element\label%
{SecProof}}

In this section we prove that the inf-sup constant for the pressure wired
Stokes element allows \W{for} a lower bound that is independent of the mesh width
and the polynomial degree $k$, and we examine the dependence on the
geometric quantity $\Theta _{\min }$ and the control parameter $\eta $. This
is formulated as the main theorem of this section.

\begin{theorem}[inf-sup stability]
\label{Theomain} 
\W{For any $k\geq 4$ and $0\leq \eta\leq 1 $,}
the inf-sup constant (\ref{infsupcond}) has the positive lower bound 
\begin{equation}
\beta \left( \mathbf{S}_{k,0}\left( \mathcal{T}\right) ,M_{\eta ,k-1}\left( 
\mathcal{T}\right) \right) \geq \left( \Theta _{\min }+\eta \right) c.
\label{infsupp_dep}
\end{equation}%
The constant $c>0$ exclusively depends on the shape-regularity of the mesh
and on $\Omega $ via a Friedrichs inequality. In particular, $c$ is
independent of the mesh width $h_{\mathcal{T}}$, the polynomial degree $k$, $%
\Theta _{\min }$, and $\eta $. 
\end{theorem}

By choosing $\eta =0$ we obtain the original Scott-Vogelius
element with a $k$-robust inf-sup constant (see \cite[Theorem 5.1]{AP_Locking}). Due to %
\eqref{press_incl}, \cite[Theorem 5.1]{AP_Locking} immediately yields that
the pressure-wired Stokes element is also inf-sup stable and the inf-sup
constant $\beta \left( \mathbf{S}_{k,0}\left( \mathcal{T}\right) ,M_{\eta
,k-1}\left( \mathcal{T}\right) \right) $ is also independent of the
polynomial degree $k$ and the mesh width $h_{\mathcal{T}}$. However, in
order to obtain the bounds as described in Theorem \ref{Theomain}, we have
to modify the arguments in \cite[Sections 4 \& 5]{AP_Locking}. The goal of
these modifications is to prove that there exists a right-inverse $\Pi
_{\eta ,k}:M_{\eta ,k-1}\left( \mathcal{T}\right) \rightarrow \mathbf{S}%
_{k,0}\left( \mathcal{T}\right) $ for the divergence operator $\operatorname{div}%
\left( \cdot \right) $. This is expressed in the following Lemma.

\begin{lemma}
\label{LemmaFinOdd}
	\W{Given any $k\geq 4$ and $0\leq\eta\leq 1$,} there exists a linear operator $\Pi_{\eta
,k}:M_{\eta,k-1}\left(  \mathcal{T}\right)  \rightarrow\mathbf{S}_{k,0}\left(
\mathcal{T}\right)  $ such that, for any $q\in M_{\eta,k-1}\left(
\mathcal{T}\right)  $,%
			\label{Thetamin_eta_mod}
\begin{align}
\operatorname*{div}\Pi_{\eta,k}  q  &  = q
\qquad\text{and}\qquad
\left\Vert \Pi_{\eta,k}q\right\Vert _{\mathbf{H}^{1}\left(
\Omega\right)  }    \leq\frac{C}{\Theta_{\min}+\eta}\left\Vert q\right\Vert
_{L^{2}\left(  \Omega\right)  }. 
\label{Thetamin_eta_modb}%
\end{align}
The constant $C > 0$ exclusively depends on the shape-regularity of the mesh. %
\end{lemma}

In order to prove this Lemma we \W{modify \cite[Lemma 4.2 and 4.5]{AP_Locking} to take into account the parameter $\eta$}; its
formulation uses the following notation. Enumerate the elements in $\mathcal{%
T}_{\mathbf{z}}$ counterclockwise by $K_{i}$ for $1\leq i\leq N_{\mathbf{z}}$
such that the edges in $\mathcal{E}_{\mathbf{z}}$ are given by $E_{i}=K_{i-1}\cap
K_{i}$, employing cyclic numbering convention, meaning $K_{0}=K_{N_{\mathbf{z}%
}}$ and $K_{N_{\mathbf{z}}+1}=K_{1}$. For $1\leq i\leq N_{\mathbf{z}}$, $%
\theta _{i}$ denotes the angle at $\mathbf{z}$ in $K_{i}$. The vectors $%
\mathbf{t}_{i}$ and $\mathbf{n}_{i}$ denote the unit tangent vector and the unit normal vector of the edge $E_{i}$ pointing away from $\mathbf{z}$ and into $K_{i}$ respectively.

\begin{lemma}[{\cite[Lemma 4.2 and 4.5]{AP_Locking}}] \label{Lem:ModAP_Lock}
	\W{Given any $k\geq 4$ and any $0\leq \eta\leq 1$, then for all}
	$\mathbf{z} \in \mathcal{V} \left(\mathcal{T}\right)$ and 
	$q \in M_{\eta,k-1} \left(\mathcal{T}\right)$, there exists a solution 
	$\left(\tilde{\mathbf{d}}_{E_i}^{\mathbf{z}}\right)_{i=1}^{N_{\mathbf{z}}} \subseteq \mathbb{R}^2$ of the system
\begin{equation} \label{eq:CompetabilityCondition}
\left. q \right\vert_{K_i} \left(\mathbf{z}\right) \sin \theta_i = \left\langle \tilde{\mathbf{d}}_{E_i}^{\mathbf{z}}, \mathbf{n}_{i+1} \right\rangle - \left\langle \tilde{\mathbf{d}}_{E_{i+1}}^{\mathbf{z}}, \mathbf{n}_i \right\rangle \qquad \forall 1 \leq i \leq N_{\mathbf{z}},
\end{equation}
that satisfies
\begin{equation} \label{eq:NormEstimate}
\sum_{i=1}^{N_{\mathbf{z}}} \left\Vert \tilde{\mathbf{d}}_{E_i}^{\mathbf{z}} \right\Vert^2 \leq C
\sum_{i=1}^{N_{\mathbf{z}}} \left. q \right\vert_{K_i} \left(\mathbf{z}\right)\W{^2} \times \begin{cases}
		\left(\Theta_{\min} + \eta\right)^{-\W{2}}, & \text{if } \; \mathbf{z} \notin \W{\mathcal{C}_{\mathcal{T}} \left(\eta\right)}, \\
1, & \text{if } \; \mathbf{z} \in \W{\mathcal{C}_{\mathcal{T}} \left(\eta\right)},
\end{cases}
\end{equation}
where $C>0$ solely depends on the shape-regularity of $\mathcal{T}$.
\end{lemma}

\begin{proof}
Let us first consider $\mathbf{z} \notin \mathcal{C}_{\mathcal{T}} \left(\eta\right)$. If $\mathbf{z} \in \mathcal{V}_{\Omega} \left(\mathcal{T}\right)$, then by \cite[Lemma 4.2]{AP_Locking} we know that there exists $\left(\mathbf{d}_{E_i}^{\mathbf{z}}\right)_{i=1}^{N_{\mathbf{z}}} \subseteq \mathbb{R}^2$ satisfying \eqref{eq:CompetabilityCondition} and
\[
\sum_{i=1}^{N_{\mathbf{z}}} \left\Vert \mathbf{d}_{E_i}^{\mathbf{z}} \right\Vert^2 \leq 2\left(1 + \frac{N_{\mathbf{z}}}{\xi \left(\mathbf{z}\right)} \right)^2 \sum_{i=1}^{N_{\mathbf{z}}} \left. q \right\vert_{K_i} \left(\mathbf{z}\right)^2 ,
\]
where $\xi \left(\mathbf{z}\right): = \sum_{i=1}^{N_{\mathbf{z}}} \left\vert \sin \left(\theta_i + \theta_{i+1}\right)
\right\vert$ employing cyclic numbering convention. We observe that from the definition of $\xi$ and \W{$\mathbf{z} \notin
\mathcal{C}_{\mathcal{T}} \left(\eta\right)$}, it follows that 
\begin{equation*}
\frac{\Theta_{\min}+ \eta}{2} \leq \Theta \left(\mathbf{z}\right) \leq \xi \left(\mathbf{z}\right).
\end{equation*}
\W{Since the maximal possible number $N_{\mathbf{z}}$ of triangles around a vertex $\mathbf{z}\in
\mathcal{V}(\mathcal{T})$ only depends on the shape-regularity,} we conclude that
\[
\sum_{i=1}^{N_{\mathbf{z}}} \left\Vert \mathbf{d}_{E_i}^{\mathbf{z}} \right\Vert^2 \leq \frac{C }{\left(\Theta_{\min} + \eta\right)^2} \sum_{i=1}^{N_{\mathbf{z}}} \left. q \right\vert_{K_i} \left(\mathbf{z}\right)^2,
\]
for some constant $C > 0$ solely depending on the shape-regularity of $\mathcal{T}$. If $\mathbf{z} \in \mathcal{V}_{\partial \Omega} \left(\mathcal{T}\right)$ we repeat the arguments presented above on the set of vectors $\left(\mathbf{d}_{E_i}^{\mathbf{z}}\right)_{i=1}^{N_{\mathbf{z}}}$ given by \cite[Lemma 4.5 (1)]{AP_Locking} and therefore by setting $\tilde{\mathbf{d}}_{E_i}^{\mathbf{z}} := \mathbf{d}_{E_i}^{\mathbf{z}}$ as described above, we have proven \eqref{eq:CompetabilityCondition} and \eqref{eq:NormEstimate} for non-$\eta$-critical vertices. Let us now consider $\mathbf{z} \in \mathcal{C}_{\mathcal{T}} \left(\eta\right)$ and we first assume $\mathbf{z} \in \mathcal{V}_{\Omega} \left(\mathcal{T}\right)$. We choose $\tilde{\mathbf{d}}_{E_i}^{\mathbf{z}} = \delta_i \mathbf{t}_i$ and some elementary computation transforms \eqref{eq:CompetabilityCondition} into
\begin{align} \label{eq:InnerCriticalAnsatz}
\left. q \right\vert_{K_i} \left(\mathbf{z}\right) & = \delta_i + \delta_{i+1} \qquad \forall 1 \leq i \leq N_{\mathbf{z}} , 
\end{align} 
employing cyclic notation convention (cf. \cite[(4.8)]{AP_Locking}). Set $ \delta_1 := 0$ and $\delta_{i+1} := \sum_{j=1}^{i} \left(-1\right)^{i-j} \left. q \right\vert_{K_j} \left(\mathbf{z}\right)$ for all $1 \leq i \leq N_{\mathbf{z}}-1$. For $i=1$ \eqref{eq:InnerCriticalAnsatz} is trivially satisfied and for $2  \leq i \leq N_{\mathbf{z}}-1$, we compute
\begin{align*}
\delta_i + \delta_{i+1} & = \sum_{j=1}^{i-1} \left(-1\right)^{i-1-j} \left. q \right\vert_{K_j} \left(\mathbf{z}\right) + \sum_{j=1}^{i} =   \left(-1\right)^{i-j} \left. q \right\vert_{K_j} \left(\mathbf{z}\right) \\
& = \sum_{j=1}^{i-1} \left(-1\right)^{i-1 -j} \left. q \right\vert_{K_j} \left(\mathbf{z}\right) + \left(-1\right)^{i-j} \left. q \right\vert_{K_j} \left(\mathbf{z}\right) + \left. q \right\vert_{K_i} \left(\mathbf{z}\right) = \left. q \right\vert_{K_i} \left(\mathbf{z}\right).
\end{align*}
For $i = N_{\mathbf{z}}$ we have by assumption that $A_{\mathcal{T}, \mathbf{z}} \left(q\right) = 0$ and conclude after rearranging the terms that $\delta_{N_{\mathbf{z}}} = \left. q \right\vert_{K_{N_{\mathbf{z}}}} \left(\mathbf{z}\right)$. Since the $\tilde{\mathbf{d}}_{E_i}^{\mathbf{z}}$ are the same as in \cite[Lemma 4.3]{AP_Locking} we obtain 
\[
\sum_{i=1}^{N_{\mathbf{z}}} \left\Vert \tilde{\mathbf{d}}_{E_i}^{\mathbf{z}} \right\Vert^2 \leq%
C \sum_{i=1}^{N_{\mathbf{z}}} \left. q \right\vert_{K_i} \left(\mathbf{z}\right)^2
\]
for some constant $C>0$ solely depending on the shape-regularity of the mesh. 
If $\mathbf{z} \in \mathcal{V}_{\partial \Omega} \left(\mathcal{T}\right)$ holds, 
the same construction as in \cite[Lemma 4.3]{AP_Locking} and $A_{\mathcal{T}, \mathbf{z}}
\left(q\right) = 0$ verify \eqref{eq:CompetabilityCondition} and \eqref{eq:NormEstimate}
also for the $\eta$-critical vertices. %
\end{proof}

\begin{proof}[Proof of Lemma \ref{LemmaFinOdd}] 
	Let $q \in M_{\eta,k-1} \left(\mathcal{T}\right)$ be given. Taking the functions $\phi_{E,k}^{\mathbf{z}}$ from \cite[Lem. B.1]{AP_Locking} we set 
\[
\Psi_{h,k} := \sum_{\mathbf{z} \in \mathcal{V}\left(\mathcal{T}\right) \atop E \in \mathcal{E}\left(\mathcal{E}\right)} \mathbf{d}_E^{\mathbf{z}} \phi_{E,k}^{\mathbf{z}},
\]
where $\mathbf{d}_E^{\mathbf{z}}$ is taken as in Lem.~\ref{Lem:ModAP_Lock}. 
Mimicking the arguments in \cite[Sec. 5]{AP_Locking} \W{(with $\Psi_{h,k}$ replaced by $\Psi_{h,p}$ therein)}
in combination \W{with} Lemma~\ref{Lem:ModAP_Lock}, we obtain that there exists a
vector field $\mathbf{u}_q \in \mathbf{S}_{k,0} \left(\mathcal{T}\right)$ satisfying $\operatorname{div} \mathbf{u}_q =
q$ and
\W{
\begin{align*}
	\left\Vert \mathbf{u}_q \right\Vert_{\mathbf{H}^{1}\left(\Omega\right)} \leq C \max\left\{1,\left(\Theta_{\min} +
	\eta\right)^{-1}\right\}/2 \left\Vert q \right\Vert_{L^{2}\left(  \Omega\right)  }\leq C \left(\Theta_{\min} +
	\eta\right)^{-1} \left\Vert q \right\Vert_{L^{2}\left(  \Omega\right)  }.
\end{align*}
with $1\leq 2(\Theta_{\min}+\eta)^{-1}$ from $\Theta_{\min},\eta\leq1$ in the last step.}
\W{The constant $C>0$ is} independent of the polynomial degree $k$, the mesh width $h_{\mathcal{T}}$, the geometric
quantity $\Theta_{\min}$ and the parameter $\eta$ and exclusively depending on the shape-regularity of the mesh\W{,} and the domain $\Omega$. 
Therefore by setting $\Pi_{\eta,k}q := \mathbf{u}_q$ we have proven our statement.
\end{proof}

\begin{proof}[Proof of Theorem \ref{Theomain}]
For the
			estimate of the inf-sup constant we compute%
			\begin{align*}
		\inf_{q\in M_{\eta,k-1}\left(  \mathcal{T}\right)  %
			\atop p \neq 0  }\sup_{\mathbf{v}\in\mathbf{S}_{k,0} \left(\mathcal{T}\right) %
			\atop \mathbf{v} \neq \mathbf{0}  }\frac{\left(  q,\operatorname*{div}\mathbf{v}\right)
			_{L^{2}\left(  \Omega\right)  }}{\left\Vert \mathbf{v}\right\Vert
			_{\mathbf{H}^{1}\left(  \Omega\right)  }\left\Vert q\right\Vert
			_{L^{2}\left(  \Omega\right)  }} & \overset{\text{Lem. } \ref{Lem:ModAP_Lock}}{\geq}\inf_{q\in M_{\eta,k-1}\left(
			\mathcal{T}\right) %
			\atop q \neq 0 }\frac{\left(
			q,\operatorname*{div}\Pi_{\eta,k}q \right)  _{L^{2}\left(  \Omega\right)}
		}{\left\Vert\Pi_{\eta,k}q \right\Vert _{\mathbf{H}^{1}\left(  \Omega%
				\right)  }\left\Vert q\right\Vert _{L^{2}\left(  \Omega\right)  }} \\
		& \geq c\left(  \Theta_{\min}+\eta\right)  ,
			\end{align*}
			where $c$ only depends on the shape-regularity of the mesh and on $\Omega$ via the Friedrichs inequality.
\end{proof}
		\W{
				Main properties of a ``good'' Stokes element are the inf-sup stability, the approximation properties of
				the velocity and pressure space, and the divergence-free condition for the velocity solution.
				We end this section by a remark on the approximation properties while the next section is devoted to the
				analysis of the smallness of the divergence.
			\begin{remark}[approximation properties]
				Since the velocity space is the standard conforming finite element space, the approximation properties
				are standard.
				For the pressure space we start with an observation for $q\in M_{\eta,k-1}(\mathcal{T})$ where $0\leq
				\eta\leq \eta_0$ with $\eta_0$ from Lemma~\ref{lem:eta_config}.
				Suppose 
				$\mathbf{z}\in\mathcal{V}_{\partial\Omega}(\mathcal{T})\cap
				\mathcal{C}_{\mathcal{T}}(\eta)$ is an $\eta$-critical boundary vertex with $N_{\mathbf{z}}$ odd (type $2$ or $4$
				in Figure~\ref{fig:vertex_patch}), $A_{\mathcal{T},\mathbf{z}}(q)=0$ reveals the
				following implication 
				\begin{align*}
					q\text{ continuous in }\mathbf{z}\Longrightarrow q(\mathbf{z})=0.
				\end{align*}
				Since the exact pressure does not vanish at theses points in general, we cannot expect a good
				approximation property of $M_{\eta,k-1}(\mathcal{T})$ in neighborhoods of such vertices that can only
				occur at corners of the domain $\Omega$ by Lemma~\ref{lem:eta_config}.
				This is a drawback for both, the original Scott-Vogelius element and the pressure-wired Stokes element.
				In \cite{BGS:PressureimprovedScottVogeliusType2024}, we present a strategy to modify the pressure space in these vertices such that standard
				approximation properties hold while the discrete inf-sup stability is preserved.
	\end{remark}}
\section{Divergence estimate of the discrete velocity\label{SecSolFree}}
			
			The Scott-Vogelius pair
			$\left(  \mathbf{S}_{k,0}\left(  \mathcal{T}\right)  ,M_{0,k-1}\left(
			\mathcal{T}\right)  \right)  $ is inf-sup stable for $k\geq4$ 
			but sensitive with respect to nearly singular vertex configurations, where $0 < \Theta \left(\mathbf{z}\right) \ll 1$ for some $\mathbf{z} \in \mathcal{V} \left(\mathcal{T}\right)$. The corresponding discrete solution is divergence free. Since our $\eta$-dependent pressure space $M_{\eta,k-1}\left( \mathcal{T}\right)$ is a proper subspace of the image of the divergence operator $\operatorname*{div}:\mathbf{S}_{k,0}(\mathcal{T})\to M_{0,k-1}\left(
			\mathcal{T}\right)$ in general, we cannot expect the discrete velocity solution 
			to be pointwise divergence-free.
			However, 
			since $\lim_{\eta \to 0} M_{\eta,k-1} \left(\mathcal{T}\right) = M_{0,k-1} \left(\mathcal{T}\right)$,
			we may expect $\left\Vert \operatorname*{div}\mathbf{u}_{\mathbf{S}}\right\Vert_{L^{2}\left(  \Omega\right)}$ to be small. The main result of this section establishes an estimate of the divergence which tends to zero as $\eta\rightarrow0$
			without any regularity assumption on the continuous Stokes problem.
		\newcommand{\Omegaeta}{\Omega(\eta)}
			\W{Consider the open subset 
				\begin{align}\label{eqn:Omega_star_def}
				\Omegaeta:=\bigcup_{\mathbf{z}\in\mathcal{C}_{\mathcal{T}}(\eta)}\operatorname*{int}(\omega_{\mathbf{z}})\subset\Omega,
				\end{align}
				where $\operatorname*{int}(\omega_{\mathbf{z}})$ denotes the interior of the closed vertex patch
				$\omega_{\mathbf{z}}$ from \eqref{nodalpatch}.
			}
			\begin{theorem}[velocity control]\label{thm:vel_control}
		Let $k\geq4$ and $0\leq\eta<\eta_{0}$ be given with $\eta_{0}>0$ as in Lemma
		\ref{lem:eta_config}. Let $(\mathbf{u}, p)$ denote the solution to \eqref{varproblemstokes} for $\mathbf{F}\in
		\mathbf{H}^{-1}(\Omega)$. Then the discrete solution $\left(
		\mathbf{u}_{\mathbf{S}},p_{M}\right)$ to
		(\ref{discrStokes}) in $\left(  \mathbf{S}_{k,0}\left(  \mathcal{T}\right)
		,M_{\eta,k-1}\left(  \mathcal{T}\right)  \right)  $ satisfies
		\begin{align*}
			\left\Vert \nabla(\mathbf{u}-\mathbf{u}_\mathbf{S})\right\Vert _{\mathbb{L}^{2}\left(
				\Omega\right)  }  &  \leq C \phantom{\eta}\min_{\substack{\mathbf{v}_\mathbf{S}\in
					\mathbf{S}_{k,0}(\mathcal{T})\\\operatorname*{div}\mathbf{v}_\mathbf{S} = 0}}\left\Vert \nabla(\mathbf{u}-\mathbf{v}_\mathbf{S})\right\Vert _{\mathbb{L}^{2}\left(
			\Omega\right)  } + \eta\phantom{^2}\hspace{-.4em}\min_{q_M\in M_{\eta,k-1}(\mathcal{T})}\|p -
		q_M\|_{L^2(\W{\Omegaeta})},\\
			\left\Vert \operatorname*{div}\mathbf{u}_\mathbf{S}\right\Vert _{{L}^{2}\left(
				\Omega\right)  }  &  \leq C\eta \min_{\substack{\mathbf{v}_\mathbf{S}\in
					\mathbf{S}_{k,0}(\mathcal{T})\\\operatorname*{div}\mathbf{v}_\mathbf{S} = 0}}\left\Vert \nabla(\mathbf{u}-\mathbf{v}_\mathbf{S})\right\Vert _{\mathbb{L}^{2}\left(
		\Omega\right)  } + \eta^2\hspace{-.4em}\min_{q_M\in M_{\eta,k-1}(\mathcal{T})}\|p - q_M\|_{L^2(\W{\Omegaeta})}.
		\end{align*}
\W{
	If $\mathcal{T}$ does \emph{not} have $\eta$-critical boundary vertices
$\mathbf{z}\in\mathcal{C}_{\mathcal{T}}(\eta)\cap \mathcal{V}_{\partial\Omega}(\mathcal{T})$ with $N_{\mathbf{z}}$ odd
(type $2$ or $4$ in Figure~\ref{fig:vertex_patch}),
then $\operatorname*{div}u_{\mathbf{S}}\equiv0$ in $\Omega\setminus\Omegaeta$.}
		The constant $C > 0$  depends solely on the shape regularity of $\mathcal{T}$ and on the domain $\Omega$. In particular $C$ is independent of the mesh width $h_{\mathcal{T}}$, the polynomial degree $k$, the geometric quantity $\Theta_{\min}$ and the control parameter $\eta$.

			\end{theorem}
			\W{A key ingredient in the proof of Theorem \ref{thm:vel_control} is the following control of the divergence.
			Define the space
			\begin{align*}
				\mathbf{S}_{\eta,k,0}(\mathcal{T}):=\left\{ \left.  \mathbf{v}_{\mathbf{S}}\in\mathbf{S}_{k,0}\left(
			\mathcal{T}\right)  \ \right\vert \ A_{\mathcal{T},\mathbf{z}}\left(  \operatorname*{div}%
			\mathbf{v}_{\mathbf{S}}\right)  =0\;\text{ for all }\mathbf{z}\in\mathcal{C}_{\mathcal{T}%
			}\left(  \eta\right)
			\right\}  .
	\end{align*} }
			\noindent For a function $v\in L^2(\Omega)$ and a subspace $U\subset
			L^2(\Omega)$
			, $v\perp U$ abbreviates
			the $L^2$ orthogonality onto $U$, i.e.,
			$\left(v,u\right)_{L^2\left(\Omega\right)} = 0$ for all $u \in U$.
			The same notation applies for vector-valued functions $\mathbf{v}\in\mathbf{L}^2(\Omega)$.
			\begin{lemma}
		[divergence control]\label{lem:div_est}
		Let $k\geq4$ and $0\leq\eta<\eta_{0}$ be given with $\eta_{0}>0$ as in Lemma
		\ref{lem:eta_config}. \W{For any $\mathbf{v}_\mathbf{S}\in \mathbf{S}_{k,0}(\mathcal{T})$ with
		$\operatorname*{div}\mathbf{v}_\mathbf{S} \perp M_{\eta,k-1}(\mathcal{T})$, there exists $C_{\mathbf{S}}\in\mathbb R$ with
		\begin{enumerate}[label=(\alph*)]
			\item $\displaystyle\operatorname*{div}\mathbf{v}_{\mathbf{S}}|_{\Omega\setminus\overline{\Omegaeta}}
				\equiv C_{\mathbf{S}}\quad\text{and}\quad
				\left\Vert \operatorname*{div}\mathbf{v}_{\mathbf{S}}-C_{\mathbf{S}}\right\Vert
			_{L^{2}\left(  \Omegaeta\right)  }\leq\sqrt{16/7}\left\Vert \operatorname*{div}\mathbf{v}_{\mathbf{S}}\right\Vert
			_{L^{2}\left(  \Omegaeta\right)  }$,
			\item $\displaystyle
	\left\Vert \operatorname*{div}\mathbf{v}_{\mathbf{S}}\right\Vert
			_{L^{2}\left(  \Omega\right)  }
				\leq\left\Vert \operatorname*{div}\mathbf{v}_{\mathbf{S}}-C_{\mathbf{S}}\right\Vert
			_{L^{2}\left(  \Omegaeta\right)  }\leq C_{\operatorname*{div}}\eta\min
				_{\mathbf{w}_{\mathbf{S}}\in\mathbf{S}_{\eta,k,0}\left(  \mathcal{T}\right) }\left\Vert \nabla\left(
		\mathbf{v}_{\mathbf{S}}-\mathbf{w}_{\mathbf{S}}\right)  \right\Vert
				_{\mathbb{L}^{2}\left(  \Omegaeta\right)  }$.
	
		\end{enumerate}
	If $\mathcal{T}$ does \emph{not} contain $\eta$-critical boundary vertices
$\mathbf{z}\in\mathcal{C}_{\mathcal{T}}(\eta)\cap \mathcal{V}_{\partial\Omega}(\mathcal{T})$ with $N_{\mathbf{z}}$ odd,
then $C_{\mathbf{S}}=0$.
}
		The constant $C_{\operatorname*{div}}>0$ exclusively depends on the shape
		regularity of $\mathcal{T}$.
			\end{lemma}

			\begin{proof}[Proof of Theorem \ref{thm:vel_control}]
				The discrete formulation \eqref{discrStokes} imposes that the discrete velocity
				$\mathbf{u}_{\mathbf{S}}$ satisfies
				$\operatorname*{div}\mathbf{u}_\mathbf{S}\perp
				M_{\eta,k-1}(\mathcal{T})$.
		\W{Consider an arbitrary $\mathbf{v}_\mathbf{S}\in \mathbf{S}_{k,0}(\mathcal{T})$ with
	$\operatorname*{div}\mathbf{v}_\mathbf{S}=0$ and observe $\mathbf{v}_\mathbf{S}\in\mathbf{S}_{\eta,k,0}(\mathcal{T})$ by definition.
The}
		weak formulation \eqref{varproblemstokes} and the discrete formulation \eqref{discrStokes} for the test function
		$\mathbf{e}_\mathbf{S}:=\mathbf{u}_{\mathbf{S}}-\mathbf{v}_{\mathbf{S}}
		\in \mathbf{S}_{k,0}(\mathcal{T})$
		satisfy
		\begin{align*}
			a(\mathbf{u}-\mathbf{u}_\mathbf{S}, \mathbf{e}_\mathbf{S}) = b(\mathbf{e}_\mathbf{S},
			p-p_M)=b(\mathbf{u}_\mathbf{S}, p-q_M) \qquad \forall q_M \in M_{\eta,k-1}(\mathcal{T})%
		\end{align*}
 with $\operatorname*{div}\mathbf{e}_\mathbf{S}=\operatorname*{div}\mathbf{u}_\mathbf{S}\perp M_{\eta,k-1}(\mathcal{T})$ in the last
 step.
\W{Let $q_M\in M_{\eta,k-1}(\mathcal{T})$ be arbitrary. Lemma~\ref{lem:div_est} provides
	$\operatorname*{div}\mathbf{u}_{\mathbf{S}}|_{\Omega\setminus\overline{\Omegaeta}}\equiv C_{\mathbf{S}}\in \mathbb R$
	and
 $C_{\mathbf{S}}\perp p-q_M\in L^2_0(\Omega)$
 reveals
 \begin{align*}
	 b(\mathbf{u}_{\mathbf{S}},
	 p-q_{M})=\left(\operatorname*{div}\mathbf{u}_{\mathbf{S}}-C_{\mathbf{S}},p-q_M\right)_{L^2(\Omegaeta)}.
\end{align*}}
\W{The previous two displayed equalities} and a Cauchy inequality provide
		\begin{align*}
			\Vert \nabla(
		\mathbf{u}_{\mathbf{S}}&-\mathbf{v}_{\mathbf{S}})  \Vert
			_{\mathbb{L}^{2}\left(  \Omega\right)  }^2
			 = a \left(\mathbf{u}_{\mathbf{S}} - \mathbf{v}_{\mathbf{S}}, \mathbf{e}_{\mathbf{S}}\right) = a(\mathbf{u}-\mathbf v_\mathbf{S}, \mathbf{e}_\mathbf{S}) -
			a(\mathbf{u}-\mathbf{u}_\mathbf{S}, \mathbf{e}_\mathbf{S})\\
			&\leq\left\Vert \nabla\left(
			\mathbf{u}-\mathbf{v}_{\mathbf{S}}\right)  \right\Vert_{\mathbb{L}^{2}\left(  \Omega\right)  } \left\Vert \nabla\mathbf{e}_{\mathbf{S}}
			\right\Vert_{\mathbb{L}^{2}\left(  \Omega\right)  }
			+\|p-q_M\|_{L^2(\W{\Omegaeta})}\|\operatorname*{div}\mathbf{u}_\mathbf{S}\W{-C_{\mathbf{S}}}\|_{L^2(\W{\Omegaeta})}.
		\end{align*}
		Lemma~\ref{lem:div_est} reveals $\left\Vert \nabla\left(
		\mathbf{u}_{\mathbf{S}}-\mathbf{v}_{\mathbf{S}}\right)  \right\Vert
		_{\mathbb{L}^{2}\left(  \Omega\right)  }\leq\left\Vert \nabla\left(
		\mathbf{u}-\mathbf{v}_{\mathbf{S}}\right)  \right\Vert_{\mathbb{L}^{2}\left(  \Omega\right)  }
		+C_{\operatorname*{div}}\eta\|p-q_M\|_{L^2(\W{\Omegaeta})}$. This, $\left\Vert \nabla\left(
		\mathbf{u}-\mathbf{u}_{\mathbf{S}}\right)  \right\Vert
		_{\mathbb{L}^{2}\left(  \Omega\right)  }\leq
		\left\Vert \nabla\left(
		\mathbf{u}-\mathbf{v}_{\mathbf{S}}\right)  \right\Vert
		_{\mathbb{L}^{2}\left(  \Omega\right)  }+
		\left\Vert \nabla\left(
		\mathbf{u}_{\mathbf{S}}-\mathbf{v}_{\mathbf{S}}\right)  \right\Vert
		_{\mathbb{L}^{2}\left(  \Omega\right)  }$ from a 
		triangle inequality, and Lemma \ref{lem:div_est} \W{lead to}
		\begin{align*}
			\left\Vert \nabla\left(
			\mathbf{u}-\mathbf{u}_{\mathbf{S}}\right)  \right\Vert
			_{\mathbb{L}^{2}\left(  \Omega\right)  }\leq
			2\left\Vert \nabla\left(
			\mathbf{u}-\mathbf{v}_{\mathbf{S}}\right)  \right\Vert_{\mathbb{L}^{2}\left(  \Omega\right)  } +C_{\operatorname*{div}}\eta\|p-q_M\|_{L^2(\W{\Omegaeta})},\\
			\left\Vert \operatorname*{div}\mathbf{u}_{\mathbf{S}}\right\Vert
			_{L^{2}\left(  \Omega\right)  }\leq C_{\operatorname*{div}}\eta\left\Vert \nabla\left(
			\mathbf{u}-\mathbf{v}_{\mathbf{S}}\right)  \right\Vert_{\mathbb{L}^{2}\left(  \Omega\right)  } +C_{\operatorname*{div}}^2\eta^2\|p-q_M\|_{L^2(\W{\Omegaeta})}.
		\end{align*}
		Since $\mathbf{v}_\mathbf{S}$ and $q_M$ are arbitrary, this concludes the proof.
			\end{proof}
		\W{One remark precedes the preparation and proof of Lemma \ref{lem:div_est} in Subsections~\ref{sub:Analysis of
		singular vertices}--\ref{sub:Divergence discrete}}.
			\W{\begin{remark}[localized pressure-robustness]
				The presence of $\eta$-critical but not singular vertices
				implies that the discrete 
				pressure space $M_{\eta,k-1}(\mathcal{T})\subsetneq
				M_{0,k-1}(\mathcal{T})=\operatorname*{div}\mathbf{S}_{k,0}(\mathcal{T})$ is
				a strict subset of the discrete divergence characterised in
			\cite{vogelius1983right}, \cite{ScottVogelius}.
				Thus, full pressure-robustness \cite{JLM+:DivergenceConstraintMixed2017} cannot be
				expected in general and the velocity error may depend on the pressure.
				However, Theorem~\ref{thm:vel_control} guarantees that this pollution effect appears only locally around
				the $\eta$-critical vertices 
				and is at most linear in $\eta$; the velocity approximation is independent of the pressure restricted to 
				$\Omega\setminus\overline{\Omegaeta}$.
		\end{remark}}
			
			\subsection{Analysis for $\eta$-critical
		vertices\label{sub:Analysis of singular vertices}}
			
			This subsection analyses  the divergence of piecewise
			smooth and globally continuous functions on the vertex patch $\mathcal{T}%
			_{\mathbf{z}}$ of an $\eta$-critical vertex $\mathbf{z}\in C_{\mathcal{T}%
			}\left(  \eta\right)  $.
			
			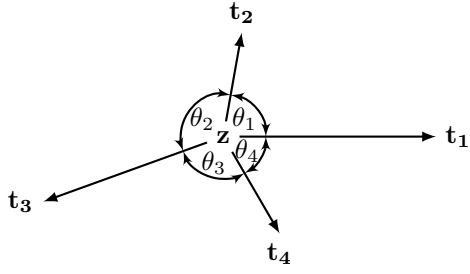
\begin{figure}[]
		\centering
		\begin{tikzpicture}[scale=0.7]
			\draw (0,0) node[] (z) {$\mathbf{z}$};
			\draw[-latex, thick] (z) -- (4,0) node (A) {};
			\node[right] at (A) {$\mathbf{t_1}$};
			\draw[-latex, thick] (z) -- ($(z)!0.5!80:(A)$) node[above] (B) {$\mathbf{t_2}$};
			\draw[-latex, thick] (z) -- ($(z)!0.9!200:(A)$) node[left] (C) {$\mathbf{t_3}$};
			\draw[-latex, thick] (z) -- ($(z)!0.53!300:(A)$) node[below] (D) {$\mathbf{t_4}$};
			\draw[latex'-latex', thick] ($(z.center)!.2!(A.center)$) arc (0:80:0.8);
			\draw[latex'-latex', thick] ($(z.center)!.2!80:(A.center)$) arc (80:200:0.8);
			\draw[latex'-latex', thick] ($(z.center)!.2!200:(A.center)$) arc (200:300:0.8);
			\draw[latex'-latex', thick] ($(z.center)!.2!300:(A.center)$) arc (300:360:0.8);
			\node at ($(40:0.8)!0.35!(z)$) {$\theta_1$};
			\node at ($(140:0.8)!0.35!(z)$) {$\theta_2$};
			\node at ($(250:0.8)!0.35!(z)$) {$\theta_3$};
			\node at ($(330:0.8)!0.35!(z)$) {$\theta_4$};
		\end{tikzpicture}
		\caption{Setting of Lemma \ref{lem:A_z_explicit}}
		\label{fig:directions}
			\end{figure}
			\begin{lemma}
		\label{lem:A_z_explicit}Consider four directions $\mathbf{t}_{1},\dots,\mathbf{t}_{4}%
		\in\mathbb{R}^{2}\backslash\left\{  \mathbf{0}\right\}  $ and some functions
		$\mathbf{v}_{1},\dots,\mathbf{v}_{4}\in\mathbf{C}^{1}\left(  U\left(
		\mathbf{z}\right)  \right)  $ defined in some neighborhood $U\left(
		\mathbf{z}\right)  $ of $\mathbf{z}\in\mathbb{R}^{2}$ as in Fig. \ref{fig:directions}. Let $\theta_{j}$
		denote the (signed) angles counted counterclockwise between $\mathbf{t}_{j}$
		and $\mathbf{t}_{j+1}$ for $j=1,\dots,4$ (with cyclic notation). If the
		G\^{a}teaux derivatives $\partial_{\mathbf{t}_{j}}(\mathbf{v}_{j-1}%
		-\mathbf{v}_{j})$
		vanish at $\mathbf{z}$ for all
		$j=1,\dots,4$ then
		\begin{equation}
			\mu\left\vert \sin\theta_{1}\right\vert \left\vert \sum_{j=1}^{4}\left(
			-1\right)  ^{j}\operatorname*{div}\mathbf{v}_{j}\left(  \mathbf{z}\right)
			\right\vert \leq\sqrt{8}\Theta\sum_{j=1}^{4}\left\Vert \nabla\mathbf{v}%
			_{j} \left(  \mathbf{z}\right)\right\Vert  \label{eqn:A_z_explicit}%
		\end{equation}
		with $\mu:=\min\left\{  \left\vert \sin\theta_{2}\right\vert ,\left\vert
		\sin\theta_{4}\right\vert \right\}  $ and $\Theta:=\max\left\{  \left\vert
		\sin\left(  \theta_{1}+\theta_{2}\right)  \right\vert ,\left\vert \sin\left(
		\theta_{2}+\theta_{3}\right)  \right\vert \right\}  $.
			\end{lemma}
			
			The proof of this lemma requires two intermediate results. Recall the
			definition of the condition number%
			\[
			\operatorname*{cond}\nolimits_{2}\left(  \mathbf{M}\right)  :=\left\Vert
			\mathbf{M}^{-1}\right\Vert _{2}\left\Vert \mathbf{M}\right\Vert _{2}%
			\]
			of a regular matrix $\mathbf{M}\in\mathbb{R}^{2\times2}$ with the induced
			Euclidean norm given by $\left\Vert \mathbf{M}\right\Vert _{2}:=\sup_{\mathbf{x}%
		\in\mathbb{R}^{2}\backslash\left\{  \mathbf{0}\right\}  }\left\Vert
			\mathbf{Mx}\right\Vert /\left\Vert \mathbf{x}\right\Vert $.
			
			\begin{proposition}
		\label{prop:cond_M}Given two vectors $\mathbf{t}_{1},\mathbf{t}_{2}%
		\in\mathbb{R}^{2}$ of unit length $\left\Vert \mathbf{t}_{1}\right\Vert
		=\left\Vert \mathbf{t}_{2}\right\Vert =1$, set $\mathbf{M}:=\W{\left[
	\mathbf{t}_{1}|\mathbf{t}_{2}\right]}  \in\mathbb{R}^{2\times2}$ (cf. Notation
		\ref{NotVectors}). Let $\theta$ be the angle between $\mathbf{t}_{1}$ and
		$\mathbf{t}_{2}$, i.e., $\cos\theta=\left\langle \mathbf{t}_{1},\mathbf{t}%
		_{2}\right\rangle $. Then
		\[
		\operatorname*{cond}\nolimits_{2}\left(  \mathbf{M}\right)  \leq\frac
		{2}{\left\vert \sin\theta\right\vert }.
		\]
		
			\end{proposition}

			\begin{proof}
			Recall the characterisation of the spectral norm $\left\Vert \mathbf M\right\Vert_{2}^{}=\sigma_{\rm max}(\mathbf M)$ as the maximal
			singular value of the matrix $\mathbf M$ from linear algebra.
			A direct computation
			reveals the eigenvalues $1\pm
			\cos\theta$ and the corresponding eigenfunctions $\left(  \pm1;1\right)\in\mathbb{R}^{2}$ of
			\[
			\mathbf{M}^{\top}\mathbf{M}=%
			\begin{pmatrix}
		1 & \cos\theta\\
		\cos\theta & 1
			\end{pmatrix}.
			\]
			Suppose $1-\cos\theta\leq1+\cos\theta$ (otherwise replace
			$\theta$ by $\pi+\theta$ and observe $\left\vert \sin\theta\right\vert
			=\left\vert \sin\left(  \theta+\pi\right)  \right\vert $). 
			Then, the estimate for the condition
			number of $\mathbf{M}$ follows from
			\[
			\operatorname*{cond}\nolimits_{2}\left(  \mathbf{M}\right)  =\sqrt
			{\frac{1+\cos\theta}{1-\cos\theta}}=\sqrt{\frac{(1+\cos\theta)^{2}}{1-\cos
			^{2}\theta}}\leq\frac{2}{|\sin(\theta)|}.
			\]\end{proof}

			The second intermediate result for the proof of Lemma \ref{lem:A_z_explicit}
			is an algebraic identity for the rotation of vectors in 2D. Define the
			rotation matrix
			\begin{equation}
		\mathbf{R}\left(  \theta\right)  :=%
		\begin{pmatrix}
			\cos\theta & -\sin\theta\\
			\sin\theta & \cos\theta
		\end{pmatrix}
		=%
		\begin{pmatrix}
			\operatorname*{Re} & -\operatorname*{Im}\\
			\operatorname*{Im} & \operatorname*{Re}%
		\end{pmatrix}
		\operatorname*{e}\nolimits^{\operatorname*{i}\theta}. \label{eqn:R_def}%
			\end{equation}

			\begin{proposition}
		[rotation identity]\label{prop:R_id} Any $\alpha,\beta,\gamma\in\mathbb{R}%
		$ satisfies%
		\[
		\sin\left(  \alpha-\beta\right)  \mathbf{R}\left(  \gamma\right)  =\sin\left(
		\gamma-\beta\right)  \mathbf{R}\left(  \alpha\right)  +\sin\left(
		\alpha-\gamma\right)  \mathbf{R}\left(  \beta\right)  .
		\]
		
			\end{proposition}

			\begin{proof}
			The real $\operatorname*{Re}$ and imaginary
			$\operatorname*{Im}$ parts are $\mathbb{R}$-linear and it
			is enough to verify
			\begin{equation}
		\sin\left(  \alpha-\beta\right)  \operatorname*{e}\nolimits^{\operatorname*{i}%
			\gamma}=\sin\left(  \gamma-\beta\right)  \operatorname*{e}%
		\nolimits^{\operatorname*{i}\alpha}+\sin(\alpha-\gamma)\operatorname*{e}%
		\nolimits^{\operatorname*{i}\beta}. \label{eqn:R_id_proof}%
			\end{equation}
			Indeed, a direct computation with the identity $2 \! \operatorname*{i}\sin
			\theta=\operatorname*{e}^{\operatorname*{i}\theta}-\operatorname*{e}%
			^{-\operatorname*{i}\theta}$ for $\theta\in\mathbb{R}$ shows
			\begin{align*}
		2 \! \operatorname*{i}\sin\left(  \gamma-\beta\right)  \operatorname*{e}%
		\nolimits^{\operatorname*{i}\alpha}  &  =\left(  \operatorname*{e}%
		\nolimits^{\operatorname*{i}\left(  \gamma-\beta\right)  }-\operatorname*{e}%
		\nolimits^{\operatorname*{i}\left(  \beta-\gamma\right)  }\right)
		\operatorname*{e}\nolimits^{\operatorname*{i}\alpha}=\operatorname*{e}%
		\nolimits^{\operatorname*{i}\gamma}\operatorname*{e}%
		\nolimits^{\operatorname*{i}\left(  \alpha-\beta\right)  }-\operatorname*{e}%
		\nolimits^{\operatorname*{i}\beta}\operatorname*{e}%
		\nolimits^{\operatorname*{i}\left(  \alpha-\gamma\right)  },\\
		2 \! \operatorname*{i}\sin\left(  \alpha-\gamma\right)  \operatorname*{e}%
		\nolimits^{\operatorname*{i}\beta}  &  =\left(  \operatorname*{e}%
		\nolimits^{\operatorname*{i}\left(  \alpha-\gamma\right)  }-\operatorname*{e}%
		\nolimits^{\operatorname*{i}\left(  \gamma-\alpha\right)  }\right)
		\operatorname*{e}\nolimits^{\operatorname*{i}\beta}=\operatorname*{e}%
		\nolimits^{\operatorname*{i}\alpha}\operatorname*{e}%
		\nolimits^{\operatorname*{i}\left(  \beta-\gamma\right)  }-\operatorname*{e}%
		\nolimits^{\operatorname*{i}\gamma}\operatorname*{e}%
		\nolimits^{\operatorname*{i}\left(  \beta-\alpha\right)  }.
			\end{align*}
			The sum of these terms and $\operatorname*{e}\nolimits^{\operatorname*{i}%
		\beta}\operatorname*{e}^{\operatorname*{i}\left(  \alpha-\gamma\right)
			}=\operatorname*{e}\nolimits^{\operatorname*{i}\alpha}\operatorname*{e}%
			^{\operatorname*{i}\left(  \beta-\gamma\right)  }$ lead to the identity
			(\ref{eqn:R_id_proof}). This and (\ref{eqn:R_def}) conclude the proof.
			\end{proof}
			\medskip
			
		\begin{proof}[Proof of Lemma \ref{lem:A_z_explicit}]
			Assume $\sin\theta_{1}\neq0$, otherwise the left-hand side of (\ref{eqn:A_z_explicit}) is zero and there is
			nothing to prove. \newline\textbf{Step 1 (preparations):} Rescale the
			directions to unit length $\left\Vert \mathbf{t}_{1}\right\Vert =\dots
			=\left\Vert \mathbf{t}_{4}\right\Vert =1$ and let $\mathbf{M}:=\W{[\mathbf{t}%
			_{1}|\mathbf{t}_{2}]}\in\mathbb{R}^{2\times2}$ be the matrix with columns
			$\mathbf{t}_{1}$ and $\mathbf{t}_{2}$. It is well known that the Piola
			transform preserves the divergence of a function. Indeed, the transformed
			function $\widehat{\mathbf{v}}:=\mathbf{M}^{-1}\mathbf{v}\circ\mathbf{\Phi}$ with
			$\mathbf{\Phi}(\mathbf{\hat{x}}):=\mathbf{M\hat{x}}+\mathbf{z}$ satisfies
			$\left(  \operatorname*{div}\mathbf{v}\right)  \circ\mathbf{\Phi}%
			=\operatorname*{div}\widehat{\mathbf{v}}$ and the componentwise G\^{a}teaux
			derivatives satisfy%
			\begin{equation}
		\partial_{\mathbf{\tilde{t}}_{j}}\left(  \widehat{\mathbf{v}}_{j-1}%
		-\widehat{\mathbf{v}}_{j}\right)  \left(  \mathbf{0}\right)  =0\quad\text{for
		}j=1,\dots,4 \label{eqn:vh_vanish}%
			\end{equation}
			in the directions $\mathbf{\tilde{t}}_{j}:=\mathbf{M}^{-1}\mathbf{t}_{j}$
			by assumption. Note
			that $\mathbf{\tilde{t}}_{1} =\mathbf{e}_{1}$ and $\mathbf{\tilde{t}}%
			_{2}=\mathbf{e}_{2}$ are the canonical unit vectors in $\mathbb{R}^{2}$ so
			that $\partial_{1}=\partial_{\mathbf{\tilde{t}}_{1}}$ and $\partial
			_{2}=\partial_{\mathbf{\tilde{t}}_{2}}$ hold. We denote the components of
			$\widehat{\mathbf{v}}_{j}$ by $\left(  \widehat{v}_{j}^{1},\widehat{v}_{j}%
			^{2}\right)  $. Following the proof of \cite[Lem. 2]{GuzmanScott2019} we
			conclude that the sum in the left-hand side of (\ref{eqn:A_z_explicit})
			equals
			\begin{align}
		\sum_{j=1}^{4}(-1)^{j}  &  \left(  \operatorname*{div}\mathbf{v}_{j}\right)
		\left(  \mathbf{z}\right)  =\sum_{j=1}^{4}(-1)^{j}\left(  \operatorname*{div}%
		\widehat{\mathbf{v}}_{j}\right)  \left(  \mathbf{0}\right)  =\sum_{j=1}%
		^{4}\left(  -1\right)  ^{j}\left(  \partial_{1}\widehat{v}_{j}^{1}%
		+\partial_{2}\widehat{v}_{j}^{2}\right)  \left(  \mathbf{0}\right) \nonumber\\
		&  =\partial_1 \! \left(\widehat{v}_{4}^{1}-\widehat{v}_{1}^{1}\right)
		\left(  \mathbf{0}\right)  +\partial_{1} \! \left(  \widehat{v}_{2}^{1}%
		-\widehat{v}_{3}^{1}\right)  \left(  \mathbf{0}\right)  +\partial_{2} \! \left(
		\widehat{v}_{2}^{2}-\widehat{v}_{1}^{2}\right)  \left(  \mathbf{0}\right)
		+\partial_{2} \! \left(  \widehat{v}_{4}^{2}-\widehat{v}_{3}^{2}\right)  \left(
		\mathbf{0}\right) \nonumber\\
		&  =\partial_{1} \! \left(  \widehat{v}_{2}^{1}-\widehat{v}_{3}^{1}\right)
		\left(  \mathbf{0}\right)  +\partial_{2} \! \left(  \widehat{v}_{4}^{2}%
		-\widehat{v}_{3}^{2}\right)  \left(  \mathbf{0}\right)  \label{eqn:A_z_bound}%
			\end{align}
			with (\ref{eqn:vh_vanish}) in the last step. The other terms might not vanish but
			(\ref{eqn:vh_vanish}) leads to a bound in terms of the full gradient
			$\nabla(\widehat{v}_{2}^{1}-\widehat{v}_{3}^{1})$ and $\nabla(\widehat{v}_{4}^{2}%
			-\widehat{v}_{3}^{2})$ at $\mathbf{0}$.\medskip\newline\textbf{Step 2 (bounds by the full
		gradient)}: The directions $\mathbf{t}_{j}=\mathbf{R}\left(  \theta_{1}%
			+\dots+\theta_{j-1}\right)  \mathbf{t}_{1}$ are rotations of $\mathbf{t}_{1}$
			by the angles $\theta_{j}$ satisfying $\cos\theta_{j}=\left\langle
			\mathbf{t}_{j},\mathbf{t}_{j+1}\right\rangle $ for $j=1,\dots,4$. The
			application of Proposition \ref{prop:R_id} for the rotation of $\mathbf{t}%
			_{1}$ first with $\alpha\leftarrow\theta_{1}+\theta_{2}$, $\beta\leftarrow0$,
			$\gamma\leftarrow\theta_{1}$ and second with $\alpha\leftarrow\theta_{1}%
			+\theta_{2}+\theta_{3}=2\pi-\theta_{4}$, $\beta\leftarrow\theta_{1}$,
			$\gamma\leftarrow0$ leads to
			\begin{align*}
		\sin\left(  \theta_{1}+\theta_{2}\right)  \,\mathbf{t}_{2}  &  =\sin\left(
		\theta_{1}\right)  \mathbf{t}_{3}+\sin\left(  \theta_{2}\right)
		\,\mathbf{t}_{1},\\
		\sin\left(  \theta_{2}+\theta_{3}\right)  \,\mathbf{t}_{1}  &  =\sin\left(
		-\theta_{1}\right)  \,\mathbf{t}_{4}+\sin\left(  2\pi-\theta_{4}\right)
		\,\mathbf{t}_{2}.
			\end{align*}
			This, (\ref{eqn:vh_vanish}) with $\mathbf{\tilde{t}}_{j}=\mathbf{e}_{j}$ for $j=1,2$, and the anti-symmetry of the
			sine result in
			\begin{equation}%
		\begin{array}
			[c]{l}%
			\left\vert \sin\left(  \theta_{2}\right)  \,\partial_{1} \! \left(  \widehat{v}%
			_{2}^{1}-\widehat{v}_{3}^{1}\right)  \right\vert =\left\vert \sin\left(
			\theta_{1}+\theta_{2}\right)  \,\partial_{2} \! \left(  \widehat{v}_{2}%
			^{1}-\widehat{v}_{3}^{1}\right)  \right\vert \leq\left\Vert \sin\left(
			\theta_{1}+\theta_{2}\right)  \,\nabla \! \left(  \widehat{v}_{2}^{1}%
			-\widehat{v}_{3}^{1}\right)  \right\Vert ,\\
			\left\vert \sin\left(  \theta_{4}\right)  \,\partial_{2} \! \left(  \widehat{v}%
			_{3}^{2}-\widehat{v}_{4}^{2}\right)  \right\vert =\left\vert \sin\left(
			\theta_{2}+\theta_{3}\right)  \,\partial_{1} \! \left(  \widehat{v}_{3}%
			^{2}-\widehat{v}_{4}^{2}\right)  \right\vert \leq\left\Vert \sin\left(
			\theta_{2}+\theta_{3}\right)  \,\nabla \! \left(  \widehat{v}_{3}^{2}%
			-\widehat{v}_{4}^{2}\right)  \right\Vert
		\end{array}
		\label{sinGateaux}%
		\end{equation}
			at $\mathbf{0}=\mathbf{\Phi}^{-1}($$\mathbf{z})$. Introduce the shorthand
			\[
			\mu:=\min\left\{  \left\vert \sin\left(  \theta_{2}\right)  \right\vert
			,\left\vert \sin\left(  \theta_{4}\right)  \right\vert \right\}
			\quad\text{and\quad}\Theta:=\max\left\{  \left\vert \sin\left(  \theta
			_{1}+\theta_{2}\right)  \right\vert ,\left\vert \sin\left(  \theta_{2}%
			+\theta_{3}\right)  \right\vert \right\}  .
			\]
			The combination of (\ref{eqn:A_z_bound}) and (\ref{sinGateaux}) provide	
			\begin{align*}
		\mu\left\vert \sum_{j=1}^{4}(-1)^{j}\left(  \operatorname*{div}\mathbf{v}%
			_{j}\right)  (\mathbf{z})\right\vert & \leq\Theta\left(  \left\Vert
			\nabla\left(  \widehat{v}_{2}^{1}-\widehat{v}_{3}^{1}\right)  \right\Vert
			+\left\Vert \nabla\left(  \widehat{v}_{3}^{2}-\widehat{v}_{4}^{2}\right)
			\right\Vert \right)  \left(  \mathbf{0}\right)  \\& \leq\sqrt{2}\Theta\sum
			_{j=1}^{4}\left\Vert \nabla\widehat{\mathbf{v}}_{j}\right\Vert \left(
			\mathbf{0}\right)
			\end{align*}
			with a triangle inequality and $\left\Vert \nabla\widehat{\mathbf{v}}_{3}%
			^{1}\right\Vert +\left\Vert \nabla\widehat{\mathbf{v}}_{3}^{2}\right\Vert
			\leq\sqrt{2}\left\Vert \nabla\widehat{\mathbf{v}}_{3}\right\Vert $ in the last step.\medskip\newline\textbf{Step 3 (finish of the
		proof):} Recall the assumption $\sin\theta_{1}\neq0$. This, the chain rule for
			derivatives $D\widehat{\mathbf{v}}_{j}=\mathbf{M}^{-1}\left(  \left(
			D\mathbf{v}_{j}\right)  \circ\mathbf{\Phi}\right)  \mathbf{M}$, and the
			Cauchy-Schwarz inequality imply $\left\Vert \nabla\widehat{\mathbf{v}}%
			_{j}\right\Vert \left(  \mathbf{0}\right)  \leq\left\Vert \mathbf{M}%
			^{-1}\right\Vert _{2}\left\Vert \mathbf{M}\right\Vert _{2}\left\Vert
			\nabla\mathbf{v}_{j}\right\Vert \left(  \mathbf{z}\right)  $ for $j=1,\dots
			,4$. Since $\operatorname*{cond}_{2}\left(  \mathbf{M}\right)  =\left\Vert
			\mathbf{M}^{-1}\right\Vert _{2}\left\Vert \mathbf{M}\right\Vert _{2}%
			\leq2/\left\vert \sin\theta_{1}\right\vert $ from Proposition \ref{prop:cond_M},
			this concludes the proof.
			\end{proof}

		Lemma \ref{lem:A_z_explicit} allows \W{for} an important generalisation of
			the known result, e.g., \cite[Lem. 2]{GuzmanScott2019}, that $A_{\mathcal{T}%
		,\mathbf{z}}\left(  \operatorname*{div}\mathbf{v}\right)  =0$ vanishes for continuous
			piecewise polynomials $\mathbf{v}\in\mathbf{S}_{k,0}\left(  \mathcal{T}%
			_{\mathbf{z}}\right)  $.
			
			\begin{corollary}
		\label{cor:A_z_explicit}
		\W{Let $k\geq0$ and $0\leq\eta<\eta_{0}$ be given with $\eta_{0}>0$ as in Lemma
		\ref{lem:eta_config}.
	Any 
		$\eta$-critical vertex $\mathbf{z}\in\mathcal{C}_{\mathcal{T}}\left(
		\eta\right)  $ satisfies}%
		\[
		\left(  \sin\phi_{\mathcal{T}}\right)  ^{2}\left\vert A_{\mathcal{T}%
			,\mathbf{z}}\left(  \operatorname*{div}\mathbf{v}\right)  \right\vert \leq C_{\operatorname*{inv}}h_{\mathbf{z}}%
		^{-1}k^{2}\eta\left\Vert \nabla\mathbf{v}\right\Vert _{\mathbb{L}^{2}\left(
			\omega_{\mathbf{z}}\right)  }\text{ for all }\mathbf{v}\in\mathbf{S}%
		_{k,0}\left(  \mathcal{T}\right)  .
		\]
		The constant $C_{\operatorname*{inv}}>0$ exclusively depends on the shape-regularity.
			\end{corollary}

			\begin{proof}
			For $\mathbf{z}\in\mathcal{V}\left(  \mathcal{T}%
			\right)  $, recall the fixed counterclockwise numbering \eqref{eqn:enum_Tz} of the triangles in
			$\mathcal{T}_{\mathbf{z}}$ so that $\mathcal{T}_{\mathbf{z}}=\left\{
			K_{j}:1\leq j\leq N_{\mathbf{z}}\right\}  $ for $N_{\mathbf{z}}:=\left\vert
			\mathcal{T}_{\mathbf{z}}\right\vert $ as in Fig. \ref{fig:vertex_patch}.
			First we consider the case of an inner 
			$\eta$-critical vertex $\mathbf{z}\in\mathcal{V}_{\Omega}\left(  \mathcal{T}%
			\right)  $ \W{with $N_{\mathbf{z}}=4$ from Lemma~\ref{lem:eta_config}(a)}. Since $\mathbf{v}\in \mathbf{S}_{k,0}\left(\mathcal{T}\right)  $ is globally continuous, the jump $\left[  \mathbf{v}\right]
			_{E_{j}}=0$ vanishes along the common edge $E_{j}:=\partial K_{j-1}%
			\cap\partial K_{j}$ and, as a consequence, $\partial_{\mathbf{t}_{j}}(
			\left.  \mathbf{v}\right\vert _{K_{j-1}}-\left.  \mathbf{v}\right\vert
			_{K_{j}})  \left(  \mathbf{z}\right)  =\mathbf{0}$. Thus, Lemma
			\ref{lem:A_z_explicit} applies to $\mathbf{v}_{j}:=\left.  \mathbf{v}%
			\right\vert _{K_{j}}$ for $j=1,\dots,4$ and shows
			\[
			\mu\left(  \sin\theta_{1}\right)  \left\vert A_{\mathcal{T},\mathbf{z}}\left(
			\operatorname*{div}\mathbf{v}\right)  \right\vert \leq\sqrt{8}\Theta\sum
			_{K\in\mathcal{T}_{\mathbf{z}}}\left\Vert \nabla \mathbf{v}%
			\vert _{K}  \left(  \mathbf{z}\right)\right\Vert .
			\]
			The $k$-explicit inverse inequality \cite[Thm.~4.76]{SchwabhpBook} and a
			scaling argument imply that there is a constant $\tilde{C}%
			_{\operatorname*{inv}}>0$ exclusively depending on the shape-regularity of
			$\mathcal{T}$ with
			\[
			\left\Vert \nabla\mathbf{v}\vert _{K}\left(
			\mathbf{z}\right)  \right\Vert \leq\left\Vert \nabla\mathbf{v}\right\Vert _{\mathbb{L}%
		^{\infty}\left(  K\right)  }\leq\tilde{C}_{\operatorname*{inv}}h_{\mathbf{z}%
			}^{-1}k^{2}\left\Vert \nabla\mathbf{v}\right\Vert _{\mathbb{L}^{2}\left(
		K\right)  } \quad \forall K\in\mathcal{T}_{\mathbf{z}}.
			\]
			Since $\Theta($$\mathbf{z})=\Theta$ holds and the minimal angle property
			implies $\sin\phi_{\mathcal{T}}\leq\sin\theta_{1}$, this, the definition of $\mu$ in
			Lemma \ref{lem:A_z_explicit}, and $\sum_{K\in\mathcal{T}_{\mathbf{z}}%
			}\left\Vert \nabla\mathbf{v}\right\Vert _{\mathbb{L}^{2}(K)}\leq2\left\Vert
			\nabla\mathbf{v}\right\Vert _{\mathbb{L}^{2}\left(  \omega_{\mathbf{z}%
		}\right)  }$ from a Cauchy inequality conclude the proof with
			$C_{\operatorname*{inv}}:=\sqrt{32}\tilde{C}_{\operatorname*{inv}}$. The case of a boundary vertex $\mathbf{z}\in\mathcal{V}_{\partial\Omega
		}\left(  \mathcal{T}\right)  $ \W{with $N_{\mathbf{z}}\leq 3$ from Lemma~\ref{lem:eta_config}(b)} can be transformed to the first case. One
			extends the \textquotedblleft open\textquotedblright\ boundary patch
			$\mathcal{T}_{\mathbf{z}}$ to a closed patch by defining shape regular
			triangles $K_{j}$, $N+1\leq j\leq4$, such that the extended patch
			$\widetilde{\mathcal{T}}_{\mathbf{z}}:=\mathcal{T}_{\mathbf{z}}\cup\left\{
			K_{j}:N+1\leq j\leq4\right\}  $ satisfies: a) $\widetilde{\mathcal{T}%
			}_{\mathbf{z}}$ a closed patch, i.e., $\mathbf{z}$ is a vertex of $K_{j}$,
			$1\leq j\leq4$ and the intersection $K_{j-1}\cap K_{j}$ is a common edge, b)
			$\mathbf{z}$ is a $\eta$-critical vertex in $\widetilde{\mathcal{T}%
			}_{\mathbf{z}}$. The function $\mathbf{v}\in\mathbf{S}_{k,0}\left(
			\mathcal{T}\right)  $ then extends to the full patch $%
			{\textstyle\bigcup\nolimits_{N+1\leq j\leq4}}
			K_{j}$ by $\mathbf{0}$ and the proof of the first case carries over.
			\end{proof}

			\subsection{Proof of Lemma \ref{lem:div_est}\label{sub:Divergence discrete}}
			
			The last ingredient for the proof of Lemma \ref{lem:div_est} concerns the
			explicit characterisation of the functions in $M_{0,k-1}(\mathcal{T})$ that
			are $L^{2}$ orthogonal to $M_{\eta,k-1}\left(  \mathcal{T}\right)  $.
			Recall the fixed counter-clockwise numbering \eqref{eqn:enum_Tz} of $\mathcal{T}_\mathbf{z}$ and let
			$\lambda_{K_j,\mathbf{z}}$ denote the barycentric coordinate associated to $\mathbf{z}$ on $K_j\in\mathcal
			T_\mathbf{z}$.

			\begin{definition}\label{def:b_z}
		For $\mathbf{z}\in\mathcal{V}\left(  \mathcal{T}\right)  $, the function $b_{k,\mathbf{z}%
		}\in\mathbb{P}_{k}\left(  \mathcal{T}_{\mathbf{z}}\right)  $ is given by%
		\begin{equation} 
			b_{k,\mathbf{z}}:=\sum_{j=1}^{N_{\mathbf{z}}}\frac{\left(  -1\right)  ^{j\W{+k}}%
			}{\left\vert K_{j}\right\vert }P_{k}^{\left(  0,2\right)  }\left(
			1-2\lambda_{K_{j},\mathbf{z}}\right)  \chi_{K_{j}}. \label{Defbkz}%
		\end{equation}
		
			\end{definition}
			\W{The polynomials $P_{k}^{\left(  0,2\right)  }\left(
			1-2\lambda_{K,\mathbf{z}}\right)\in \mathbb P_k(K)$ 
have been
introduced in \cite{BS:GaussLegendreElementsStable2007} to characterize the orthogonal complement
of $\operatorname*{div}(\mathbf{S}_{k+1,0}(\mathcal{T}))$ in $\mathbb P_{k,0}(\mathcal{T})$ on certain macro elements, see also
\cite{CCSS_CR_1}, \cite{SauterCR_prob}, \cite{SauterTorres_CR3D}.}
			\begin{lemma}
		\label{lem:b_z_def}
		Set $\zeta_{k}:=\binom{k+2}{2}$ for $k\geq 1$ and let $\mathbf{z}\in\mathcal{V}\left(  \mathcal{T}\right)
		$.

		\begin{enumerate}
			\item The function $b_{k,\mathbf{z}}\in\mathbb{P}_{k}\left(  \mathcal{T}%
			_{\mathbf{z}}\right)  $ from (\ref{Defbkz}) satisfies, for all $1\leq j\leq
			N_{\mathbf{z}}$, that%
			\begin{align}
			\left.  b_{k,\mathbf{z}}\right\vert _{K_{j}}\left(  \mathbf{y}\right)
			=\frac{\left(  -1\right)  ^{j}}{\left\vert K_{j}\right\vert }\times\left\{
			\begin{array}
				[c]{ll}%
				\zeta_{k} &\text{if } \mathbf{y}=\mathbf{z},\\
				\W{(-1)^k} &\text{else}
			\end{array}
		\right.&&\forall\mathbf{y}\in\mathcal{V}(K_j).\label{eqn:bkz_y}
			\end{align}
			
		\item The following integral relations hold \W{on the triangle $K_j\in \mathcal{T}_{\mathbf{z}}$}:%
			\W{\begin{align}
				\left(  b_{k,\mathbf{z}},1\right)  _{L^{2}\left(  K_{j}\right)  }  &  =\left(
				-1\right)  ^{j}\zeta_{k}^{-1}, &  &  \label{bz1}\\%
				\left(  b_{k,\mathbf{z}},q\right)  _{L^{2}\left(  K_{j}\right)  }  &
				=q\left(  \mathbf{z}\right)  \left(  b_{k,\mathbf{z}},1\right)  _{L^{2}\left(
					K_{j}\right)  } & \forall q  &  \in\mathbb{P}_{k}(K_{j}).\label{bz0}
\end{align}}
\item \W{The function $\zeta_{k}b_{k,\mathbf{z}}$ is the Riesz representation of $A_{\mathcal{T},\mathbf{z}}$ in
		$\mathbb P_k(\mathcal{T}_{\mathbf{z}})$, i.e.,
\begin{align}
		\zeta_{k}\left(  b_{k,\mathbf{z}},q\right)  _{L^{2}\left(  \omega_{\mathbf{z}}\right)
		}  &  =%
		A_{\mathcal{T}%
		,\mathbf{z}}\left(  q\right)& \forall q&\in \mathbb P_{k}(\mathcal{T}_{\mathbf{z}}). \label{Eq:Value scalarproduct bkz q}
\end{align}
		In partiular,} %
		$b_{k,\mathbf{z}}$ is
			$L^{2}(  \omega_{\mathbf{z}})  $ orthogonal to all 
			$q\!\in\!\mathbb{P}_{k}(\mathcal{T}_{\mathbf{z}})$ with $A_{\mathcal{T}%
				,\mathbf{z}}(  q)  \!=\!0$.
			\item For $k\geq2$, the set $\{1\}\cup\{
				b_{k,\mathbf{z}}\}_{\mathbf{z}\in\mathcal{V}(\mathcal{T})}$ is linearly independent and
				\begin{align}\label{bk_min}
				\W{\frac{3}{4}\Bigg\|\sum^{}_{\mathbf{z}\in\mathcal{V}(K)}
					c_\mathbf{z}b_{k,\mathbf{z}} 
				\Bigg\|_{L^2(K)}^2}\hspace*{-1em}\leq|K|^{-1}\sum^{}_{\mathbf{z}\in\mathcal{V}(K)} c_\mathbf{z}^2\leq\frac{12}{7}
					\operatorname*{min}_{C\in\mathbb R}\Bigg\|\sum^{}_{\mathbf{z}\in\mathcal{V}(K)}
					c_\mathbf{z}b_{k,\mathbf{z}} -
					C\Bigg\|_{L^2(K)}^2
				\end{align}
			\W{for any $K\in\mathcal{T}$ and $c_\mathbf{z}\in\mathbb R$.}
		\end{enumerate}
			\end{lemma}
			\begin{proof}
			From \cite[(3.14), (3.2)]{CCSS_CR_1} it follows that the function
			$b_{j}:=P_{k}^{\left(  0,2\right)  }\left(  1-2\lambda_{K_{j},\mathbf{z}%
			}\right)  $ is orthogonal to any function $q\in\mathbb{P}_{k}(K_{j})$ with $q\left(  \mathbf{z}\right)  =0$ and fulfils
			\begin{equation}
		b_{j}\left(  \mathbf{z}\right)  =P_{k}^{(0,2)}\left(  -1\right)  =\left(
		-1\right)  ^{k}\zeta_{k},\quad b_{j}\left(  \mathbf{y}\right)
		=P_{k}^{\left(  0,2\right)  }\left(  1\right)  =1\quad\forall\mathbf{y}%
		\in\mathcal{V}\left(  K_{j}\right)  \backslash\left\{  \mathbf{z}\right\}  .
		\label{bKvertexvalues}%
			\end{equation}
			The integral of $b_{j}$ over the triangle $K_j$ can be evaluated explicitly as%
			\begin{align}
		\left(  b_{j},1\right)  _{L^{2}\left(  K_{j}\right)  }&
		  =\int_{K_{j}}P_{k}^{\left(  0,2\right)  }\left(
		1-2\lambda_{K_{j},\mathbf{z}}\right)  =2\left\vert
		K_{j}\right\vert \int_{0}^{1}\int_{0}^{1-x_{1}}P_{k}^{\left(  0,2\right)
		}\left(  1-2x_{1}\right)  d x_{2} d x_{1}\nonumber\\
		&  =2\left\vert K_{j}\right\vert \int_{0}^{1}\!\!\left(
		1-x_{1}\right)  P_{k}^{\left(  0,2\right)  }(  1-2x_{1})
	d x_{1}
		 \nonumber\\ & \overset{t=1-2x_{1}}{=}\frac{|K_j|}{2} \int_{-1}^{1}\!\!\left(  t+1\right)  P_{k}^{\left(  0,2\right)
		}(  t)  d t\nonumber\\
		&  \overset{\text{\cite[Lem. C.1]{CCSS_CR_1}}}{=}\left(  -1\right)  ^{k}\left\vert K_{j}\right\vert
		/\zeta_{k}. \label{bK_1_sc_p}%
			\end{align}
			This implies (\ref{bz1}).
			The $L^{2}$ orthogonality $\left(
			b_{j},q-q\left(  \mathbf{z}\right)  \right)  _{L^{2}\left(  K_{j}\right)  }=0$ %
			 shows
			\[
			0<\left(  b_{j},b_{j}\right)  _{L^{2}\left(  K_{j}\right)  }=b_{j}\left(
			\mathbf{z}\right)  \;\left(  b_{j},1\right)  _{L^{2}\left(  K_{j}\right)
			}\overset{\text{(\ref{bKvertexvalues}), (\ref{bK_1_sc_p})}}{=}\left\vert
			K_{j}\right\vert \qquad\forall q\in\mathbb{P}_{k}\left(  K_{j}\right)  .
			\]
			In turn, (\ref{bz0}) follows from
			\[
			\left(  b_{k,\mathbf{z}},q\right)  _{L^{2}\left(  K_{j}\right)  }=\left(
			b_{k,\mathbf{z}},q-q\left(  \mathbf{z}\right)  \right)  _{L^{2}\left(
		K_{j}\right)  }+\left(  b_{k,\mathbf{z}},q\left(  \mathbf{z}\right)  \right)
			_{L^{2}\left(  K_{j}\right)  }=q\left(  \mathbf{z}\right)  \left(
			b_{k,\mathbf{z}},1\right)  _{L^{2}\left(  K_{j}\right)  }.
			\]
			The definition of $\left.  b_{k,\mathbf{z}}\right\vert _{K_{j}}=\left(
				-1\right)  ^{j\W{+k}}\left\vert K_{j}\right\vert ^{-1}b_{j}$ in $K_{j}\in
			\mathcal{T}_{\mathbf{z}}$ reveals
			\begin{align*}
		\left(  b_{k,\mathbf{z}},q\right)  _{L^{2}\left(  \omega_{\mathbf{z}}\right)
		}  &  =\sum_{j=1}^{N_{\mathbf{z}}}\left.  q\right\vert _{K_{j}}\left(
		\mathbf{z}\right)  \frac{\left(  -1\right)  ^{j\W{+k}}}{\left\vert K_{j}\right\vert
		}\left(  b_{j},1\right)  _{L^{2}\left(  K_{j}\right)  } 
		\overset{\text{(\ref{bK_1_sc_p})}}{=}%
		\zeta_{k}^{-1}A_{\mathcal{T}%
			,\mathbf{z}}\left(  q\right)%
			\end{align*}
		for any $q\in\mathbb{P}_{k}(\mathcal{T}_{\mathbf{z}})$. This \W{is \eqref{Eq:Value scalarproduct bkz q} and} implies the
			orthogonality $\left(  b_{k,\mathbf{z}},q\right)  _{L^{2}(\omega_{\mathbf{z}%
		})} = 0$ for any $q\in\mathbb{P}_{k}\left(  \mathcal{T}_{\mathbf{z}}\right)  $
			with $A_{\mathcal{T},\mathbf{z}}\left(  q\right)  =0$.
			Given coefficients $c_\mathbf{z}\in\mathbb R$ for $\mathbf{z}\in\mathcal{V}(K),K\in\mathcal{T}$, set
			$q:=\sum^{}_{\mathbf{z}\in\mathcal{V}(K)} c_\mathbf{z}b_{k,\mathbf{z}}$.
			The minimum in \eqref{bk_min} is obtained for the integral mean $q_K=|K|^{-1}
			(q,1)_{L^2(K)}$, i.e.,
			\begin{align*}
				\operatorname*{min}_{C\in\mathbb R}\|q &- C\|_{L^2(K)}^2=
				\|q-q_K\|_{L^2(K)}^2=\|q\|_{L^2(K)}^2-\|q_K\|_{L^2(K)}^2=\sum^{}_{\mathbf{z},\mathbf{y}\in\mathcal{V}(K)}
				c_\mathbf{z}c_\mathbf{y}M_{\mathbf{z},\mathbf{y}}%
			\end{align*}
			with $M_{\mathbf{z},\mathbf{y}}:=(b_{k,\mathbf{z}},b_{k,\mathbf{y}})_{L^2(K)}-|K|^{-1}
					(b_{k,\mathbf{z}},1)_{L^2(K)}(b_{k,\mathbf{y}},1)_{L^2(K)}$.
					Since $\|b_{k,\mathbf{z}}\|_{L^2(K)}^2=|K|^{-1}$ and
				$|(b_{k,\mathbf{z}},b_{k,\mathbf{y}})_{L^2(K)}|=\zeta_k^{-1}|K|^{-1}$ from
				\eqref{eqn:bkz_y}--\eqref{bz0} for $\mathbf{z}\ne\mathbf{y}\in\mathcal{V}(K)$, \eqref{bz1} shows that
				the diagonal and offdiagonal elements of the
					symmetric diagonally dominant $3\times 3$ matrix
				$M:=(M_{\mathbf{z},\mathbf{y}})_{\mathbf{z},\mathbf{y}\in\mathcal{V}(K)}\in\mathbb R^{3\times3}$ are
				controlled by
\begin{align*}
M_{\mathbf{z},\mathbf{z}}=|K|^{-1}(1-\zeta_k^{-2})=:m&& 
|M_{\mathbf{z},\mathbf{y}}|\leq |K|^{-1}(\zeta_k^{-1}+\zeta_k^{-2})=:r
\end{align*}
$\mathbf{z}\ne\mathbf{y}\in \mathcal{V}(K)$. 
				This and the Gerschgorin circle theorem \cite[Thm.~7.2.1]{golub_matrix_1996} prove that (all and in particular) the lowest eigenvalue
				$\lambda_{\rm min}$ of $M$ is bounded \W{from} below by \linebreak
				$m-2r\geq7|K|^{-1}/12$ for $k\geq 2$.
				\W{Analog arguments control the maximal eigenvalue
					$\lambda_{\max}\leq|K|^{-1}(1+2\zeta_k^{-1})\leq |K|^{-1}4/3$ for
				$k\geq2$ of $N\in \mathbb R^{3\times3}$ with
			$N_{\mathbf{z},\mathbf{y}}:=(b_{k,\mathbf{z}},b_{k,\mathbf{y}})_{L^2(K)}$. }
				Hence, the min-max theorem implies \eqref{bk_min} and the linear independence follows from
			the support property, $\operatorname*{supp}b_{k,\W{\mathbf{z}}}=\omega_{\W{\mathbf{z}}}$.
				This
			concludes the
			proof.
			\end{proof}
			Now we have all the ingredients for the proof of the key Lemma \ref{lem:div_est} of this section.
			\medskip
			\begin{proof}[Proof of Lemma \ref{lem:div_est}]
				This proof is split into 3 
			steps.\medskip\newline
			\noindent\textbf{Step 1 (characterisation of the divergence)}: Note that 
			$\operatorname*{div}\mathbf{v}_{\mathbf{S}}\in \mathbb P_{k-1,0}\left(  \mathcal{T}%
			\right)  $ has integral mean zero from $\mathbf{v}_{\mathbf{S}}\in
			\mathbf{S}_{k,0}(\mathcal T)$. The $L^{2}$ orthogonality  $\operatorname*{div}\mathbf{v}_{\mathbf{S}}\perp
			M_{\eta,k-1}\left(  \mathcal{T}\right)$ implies
			\[
			\operatorname*{div}\mathbf{v}_{\mathbf{S}}\in M_{\eta,k-1}\left(
			\mathcal{T}\right)  ^{\perp}:=\left\{  p\in \mathbb P_{k-1}\left(  \mathcal{T}%
			\right)  \ :\ \left(  p,q\right)  _{L^{2}\left(  \Omega\right)  }=0\text{ for
		all }q\in M_{\eta,k-1}\left(  \mathcal{T}\right)  \right\}  .
			\]
			It follows from Lemma \ref{lem:b_z_def}(4) that $%
\operatorname*{span}\{1,b_{k-1,\mathbf{z}}\ :\ \mathbf{z}\in
			C_{\mathcal{T}}\left(  \eta\right) \}$ is the orthogonal compliment of $M_{\eta,k-1}\left(  \mathcal{T}\right)$ in $\mathbb{P}_{k-1} \left(\mathcal{T}\right)$. This guarantees 
			the existence of coefficients $c_{\mathbf{z}}\in\mathbb{R}$ for $\mathbf{z}\in\mathcal{C}_{\mathcal
			T}\left(\eta\right)  $ and $C_{\mathbf{S}}\in\mathbb R$ such that
			\begin{align}
				\operatorname*{div}\mathbf{v}_{\mathbf{S}}=\sum_{\mathbf{z}\in\mathcal{C}_{\mathcal T}%
				\left(  \eta\right)  }c_{\mathbf{z}} b_{k-1,\mathbf{z}} +C_{\mathbf{S}}.\label{eqn:div_char}
			\end{align}
			\W{Recall the open subset $\Omegaeta\subset\Omega$ from \eqref{eqn:Omega_star_def}.
				Since $\operatorname*{supp}b_{k-1,\mathbf{z}}=\omega_{\mathbf{z}}\subset\overline{\Omegaeta}$ from Definition~\ref{def:b_z},
				\eqref{eqn:div_char} verifies $\operatorname*{div}
				\mathbf{v}_{\mathbf{S}}|_{\Omega\setminus\overline{\Omegaeta}}\equiv C_{\mathbf{S}}$.
				For any $K\in \mathcal{T}$, Lemma~\ref{lem:b_z_def}(4) 
				provides
				\begin{align*}
					\|\operatorname*{div}\mathbf{v}_{\mathbf{S}}-C_{\mathbf{S}}\|_{L^2(K)}=\left\|\sum^{}_{\mathbf{z}\in
					\mathcal{C}_{\mathcal{T}}(\eta)} c_{\mathbf{z}}b_{k-1,\mathbf{z}}\right\|_{L^2(K)}\leq
					\sqrt{16/7}\|\operatorname*{div}\mathbf{v}_{\mathbf{S}}\|_{L^2(K)}.
				\end{align*}
				The integral mean is the best-approximation onto constants and vanishes for $\operatorname*{div}\mathbf{v}_{\mathbf{S}}$
				so that this and $\operatorname*{div}\mathbf{v}_{\mathbf{S}}|_{\Omega\setminus\overline{\Omegaeta}}\equiv C_{\mathbf{S}}$ verify
				\begin{align}\nonumber
					\|\operatorname*{div}\mathbf{v}_{\mathbf{S}}\|_{L^2(\Omega)}
					&=\min_{C\in\mathbb
					R}\|\operatorname*{div}\mathbf{v}_{\mathbf{S}}\!-\!C\|_{L^2(\Omega)}\\
					&\leq\|\operatorname*{div}\mathbf{v}_{\mathbf{S}}
					-C_{\mathbf{S}}\|_{L^2(\Omegaeta)}
					\leq\sqrt{16/7}\|\operatorname*{div}\mathbf{v}_{\mathbf{S}}\|_{L^2(\Omegaeta)}.\label{eqn:div_C_est}
				\end{align}
				Lemma~\ref{lem:eta_config}(a) and $\eta\leq \eta_0$ guarantee $N_{\mathbf{z}}=4$ for any interior $\eta$-critical
			vertex $\mathbf{z}\in \mathcal{C}_{\mathcal{T}}(\eta)\cap\mathcal{V}_{\Omega}(\mathcal{T})$.
		If $\mathcal{T}$ does not contain $\eta$-critical boundary vertices $\mathbf{z}\in
		\mathcal{C}_{\mathcal{T}}(\eta)\cap\mathcal{V}_{\partial\Omega}(\mathcal{T})$ with $N_{\mathbf{z}}$ odd, the integral mean
	of $b_{k-1,\mathbf{z}}$ vanishes for all $\mathbf{z}\in\mathcal{C}_{\mathcal{T}}(\eta)$ by
	\eqref{bz1} and 
\begin{align*}
	|\Omega|\,C_{\mathbf{S}}=\int_\Omega \operatorname*{div}\mathbf{v}_{\mathbf{S}}-
	\sum^{}_{\mathbf{z}\in\mathcal{C}_{\mathcal{T}}(\eta)} c_{\mathbf{z}}\int_\Omega b_{k-1,\mathbf{z}} = 0.
\end{align*}}
			\textbf{Step 2 (preparations)}: 
			\W{Let $K_{\mathbf{z}}\in
			\mathcal{T}_{\mathbf{z}}$ denote any fixed triangle in the vertex patch of $\mathbf{z}\in \mathcal{C}_{\mathcal{T}}(\eta)$.}
		\W{The orthogonality of $\operatorname*{div}\mathbf v_\mathbf{S}$ onto $C_{\mathbf{S}}\in\mathbb R$,
	\eqref{eqn:div_char}, and Lemma~\ref{lem:b_z_def}(3) provide}
			\begin{align}
				\left\Vert \operatorname*{div}\mathbf{v}_{\mathbf{S}}\right\Vert
			_{L^{2}\left(  \Omega\right)  }^{2}   & =\sum_{\mathbf{z}\in C_{\mathcal{T}}(  \eta)  }\left(  c_{\mathbf{z}}b_{k-1,\mathbf{z}%
		},\operatorname*{div}\mathbf{v}_{\mathbf{S}}\right)  _{L^{2}\left(
			\omega_{\mathbf{z}}\right)  }\W{\overset{\eqref{Eq:Value scalarproduct bkz q}}{=}{\zeta_{k-1}^{-1}}\sum_{\mathbf{z}\in C_{\mathcal{T}}(  \eta)
		}c_{\mathbf{z}}\,A_{\mathcal{T},\mathbf{z}}\left(  \operatorname*{div}%
			\mathbf{v}_{\mathbf{S}}\right)}\nonumber\\
			& \W{\leq\zeta_{k-1}^{-1}\sqrt{\sum_{\mathbf{z}\in \mathcal C_{\mathcal{T}}\left(  \eta\right)
			}\left\vert K_{\mathbf{z}}\right\vert %
				\left(  A_{\mathcal{T},\mathbf{z}}\left(
			\operatorname*{div}\mathbf{v}_{\mathbf{S}}\right)  \right)  ^{2}}\sqrt
		{\sum_{\mathbf{z}\in C_\mathcal{T}(\eta)  }\left\vert K_{\mathbf{z}}\right\vert ^{-1}c_\mathbf{z}^2}}%
		\label{eqn:proof_divv_h}%
			\end{align}
			\W{with a Cauchy inequality in $\ell^2$ %
			in the last step. Lemma~\ref{lem:b_z_def}(4) controls the last term by
			\begin{align}\label{eqn:div_bound}
			\sum_{\mathbf{z}\in C_\mathcal{T}(\eta)  }\left\vert K_{\mathbf{z}}\right\vert ^{-1}c_\mathbf{z}^2\leq
		\sum^{}_{K\in \mathcal{T}} \sum^{}_{\mathbf{z}\in \mathcal{V}(K)\cap\mathcal{C}_{\mathcal{T}}(\eta)}
			|K|^{-1}c_{\mathbf{z}}^2\leq12/7\|\operatorname*{div}\mathbf{v}_{\mathbf{S}}\|_{L^2(\Omegaeta)}^2.
	\end{align}
}
	\begin{figure}%
		\centering
		\includegraphics[width=0.3\textwidth]{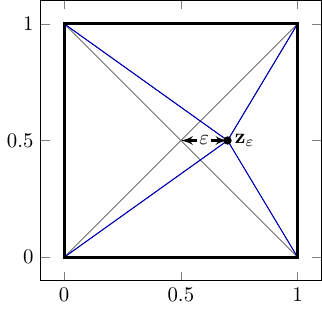}
		\caption{Perturbation $\mathcal T_\varepsilon$ (blue) of the criss-cross triangulation $\mathcal{T}_0$ (gray) of the square}
		\label{fig:Square_init_eps}
			\end{figure}			
			\noindent\textbf{Step 3 (\W{divergence control})}: 
		Let $\mathbf w_\mathbf{S}\in \W{\mathbf S_{\eta,k,0}}(\mathcal T)\W{\subset\mathbf{S}_{k,0}(\mathcal{T})}$ be arbitrary and recall $\left\vert K\right\vert
			\leq h_{\mathbf{z}}^{2}$ for all $K\in\mathcal{T}_\mathbf{z}$ by definition. 
			Since $A_{\mathcal{T},\mathbf{z}}\left(
			\operatorname*{div}\mathbf{w}_{\mathbf{S}}\right)$ vanishes and $\Theta\left(  \mathbf{z}\right)  \leq\eta$ for every
			$\mathbf{z}\in \mathcal C_{\mathcal{T}}\left(  \eta\right)$,%
			(\ref{eqn:proof_divv_h}%
			)--(\ref{eqn:div_bound}) and Corollary~\ref{cor:A_z_explicit} applied to
			$\mathbf{v}_{\mathbf{S}}-\mathbf{w}_{\mathbf{S}}\in\mathbf{S}_{k,0}\left(
			\mathcal{T}\right)  $ imply
			\begin{align*}
		\left\Vert \operatorname*{div}\mathbf{v}_{\mathbf{S}}\right\Vert
		_{L^{2}\left(  \Omega\right)  }  &  \leq \sqrt{12/7}\zeta_{k-1}^{-1}\sqrt{\sum_{\mathbf{z}\in \mathcal C_{\mathcal{T}%
				}\left(  \eta\right)
			}h_{\mathbf{z}}^{2}\left(  A_{\mathcal{T},\mathbf{z}}\left(
			\operatorname*{div}\mathbf{v}_{\mathbf{S}}-\operatorname*{div}\mathbf{w}%
			_{\mathbf{S}}\right)  \right)  ^{2}}\\
		&  \leq \sqrt{12/7}\zeta_{k-1}^{-1}k^{2}%
		\sqrt{\sum_{\mathbf{z}\in \mathcal C_{\mathcal{T}}\left(  \eta\right)}\Theta\left(  \mathbf{z}\right)
			^{2}\left\Vert \nabla\left(  \mathbf{u}_{\mathbf{S}}-\mathbf{w}_{\mathbf{S}%
			}\right)  \right\Vert _{\mathbb{L}^{2}\left(  \omega_{\mathbf{z}}\right)
			}^{2}}.
			\end{align*}
			\W{Hence,
				the finite overlap of the
			vertex patches results in
			\begin{align*}
				\|\operatorname*{div}\mathbf{v}_{\mathbf{S}}\|_{L^2(\Omega)}
				\leq \frac{C_{\operatorname*{inv}}}{\sin^{2}\phi_{\mathcal{T}}%
		}\sqrt{36/7}\zeta_{k-1}^{-1}k^{2}\eta\left\Vert \nabla\left(  \mathbf{u}_{\mathbf{S}}-\mathbf{w}%
		_{\mathbf{S}}\right)  \right\Vert _{\mathbb{L}^{2}\left(  \Omegaeta\right)}.
\end{align*}}
Since $\mathbf{w}_{\mathbf{S}}\in\W{\mathbf S_{\eta,k,0}}(\mathcal T)$ was arbitrary, \W{this, \eqref{eqn:div_C_est}, and
		$\zeta_{k-1}^{-1}k^2=2k/(k+1)\leq 2$}
			conclude the proof.
			\end{proof}
			\section{Numerical experiments\label{sec:Numerical experiments}}
			
			In this section we report on numerical experiments of the convergence rates
			for the pressure-wired Stokes elements and investigate the dependence of $\left\Vert \operatorname*{div}\mathbf{u}_{\mathbf{S}}\right\Vert _{L^{2}\left( \Omega\right)  }$ on $\eta$ and the number of degrees of freedom.

\begin{figure}%
		\centering
		\hbox{\hspace{-4em}\includegraphics[width=1.3\textwidth]{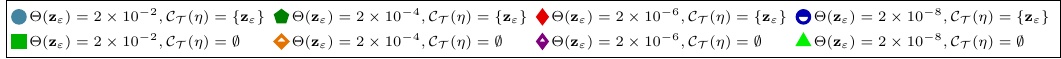}}
		\hbox{\hspace{-2em}\includegraphics{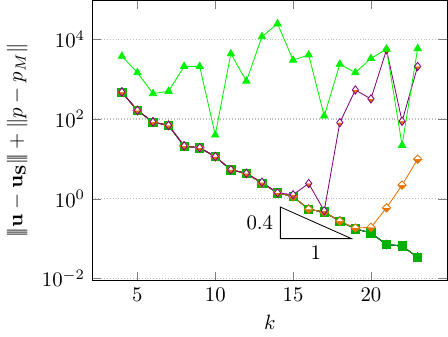}\hspace{1em}
			\includegraphics{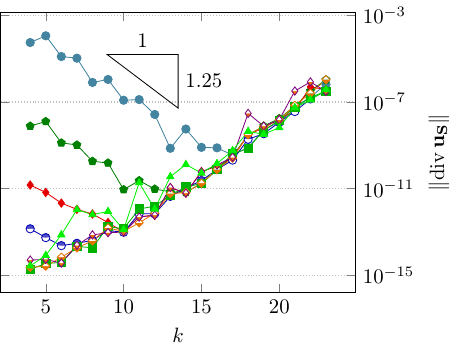}}
		\caption{Convergence history of the discretisation error (left) and $\|\operatorname{div}\,
			\mathbf{u}_{\mathbf{S}}\|$ (right) with $\mathbf{z}_\varepsilon$ treated as nearly singular
			$\mathbf{z}_\varepsilon\in\mathcal{C}_{\mathcal{T}}(\eta)$ or non-singular
			$\mathbf{z}_\varepsilon\not\in\mathcal{C}_{\mathcal{T}}(\eta)$ vertex for the $k$-method}
		\label{fig:Square_exact_p}
			\end{figure}

\W{We comment on a straightforward implementation with a direct LU solver and Lagrange multipliers for the pressure
	constraints.
For the original Scott-Vogelius element and
nearly singular vertices the condition number of the algebraic system
becomes very large.
As a consequence all input errors and variational crimes
are amplified by the bad conditioning.
The treatment of the ill-conditioning in the presence of few nearly singular vertices by more robust solvers, e.g, of Krylov
type, lies beyond the scope of this
paper with its focus
on the pressure-wired Stokes element.}
			\newline\textbf{Analytical solution on the unit square.}
			The difference to the classical Scott-Vogelius element stems from the different treatment of near-singular
			vertices only.
			Therefore we focus our benchmark on the criss-cross triangulation $\mathcal{T}_{\varepsilon}$ of the unit square
			$\Omega:=(0,1)^2$ with interior vertex $\mathbf{z}_{\varepsilon}:=(1/2+\varepsilon;1/2)$ perturbed by some $0<\varepsilon<1/2$ shown in Fig. \ref{fig:Square_init_eps}.
			For $\varepsilon\ll1$, this triangulation locally models a critical mesh configuration with
			$0<\Theta(\mathbf{z}_\varepsilon)\ll1$ where in finite
			arithmetic the classical
			Scott-Vogelius FEM becomes unstable.
			Consider the exact smooth solution
			\begin{align*}
			\mathbf{u}\left(  \mathbf{x}\right)
			=\begin{pmatrix}
		\phantom{-}\sin^{2}\left(  \pi x_{1}\right)  \sin\left(  \pi x_{2}\right)\cos\left(  \pi x_{2}\right)\\
		-\sin^{2}\left(  \pi x_{2}\right)  \sin\left(  \pi x_{1}\right)\cos\left(  \pi x_{1}\right)
	\end{pmatrix} ,&& \hspace{-.5em}p\left(
			\mathbf{x}\right)  =10^{6}\operatorname*{e}\nolimits^{-\left(  x_{1}%
		-0.3\right)  ^{-2}-\left(  x_{2}-32/500\right)  ^{-2}}\hspace{-.4em}+C,
			\end{align*}
			to the Stokes problem with body force $\mathbf{F}\equiv\mathbf{f}:=-\Delta\mathbf{u}%
			+\nabla p\in\mathbf{C}^{\infty}\left(  \Omega\right)  $.
			The velocity is pointwise divergence-free as it is the vector curl of $\sin^{2}\left(  \pi x_{1}\right)  \sin^2\left(
			\pi x_{2}\right)$ and $C$ is chosen such that the pressure $p\in L^2_0(\Omega)$ has zero integral mean.
			\newline\textbf{Optimal convergence of the $h$-method.}
			This benchmark considers uniform red-refinement of the initial mesh $\mathcal T_\varepsilon$ that subdivides each
			element into four congruent children by joining the edge midpoints.
			This refinement strategy does not introduce new near-singular vertices in the refinement and
			$\Theta(\mathbf{z}_\varepsilon)$ remains constant throughout the refinement.
			We consider $\varepsilon=0.01,0.0001,10^{-6},10^{-8}$ that corresponds to $\Theta(\mathbf{z}_{\varepsilon})
			\approx 2\varepsilon$.
			\begin{figure}[ptb]
		\centering
	\includegraphics[width=\textwidth]{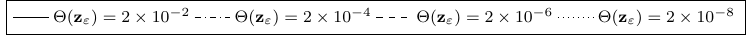}\\
	\hbox{\hspace*{-2em}\includegraphics[width=1.2\textwidth]{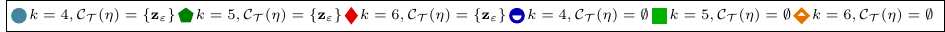}}
		\hbox{\hspace*{-2em}\includegraphics{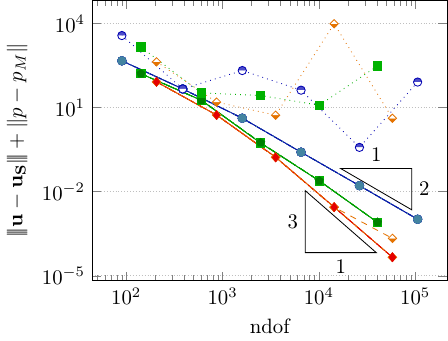}\hspace{1em}
			\includegraphics{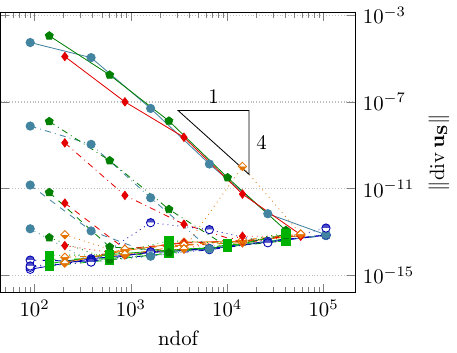}}
		\caption{Convergence history of the total error (left) and of the divergence (right) with
			$\mathbf{z}_\varepsilon$ treated as \W{$\eta$-critical} 
			$\mathbf{z}_\varepsilon\in\mathcal C_\mathcal{T}(\eta)$ or non-singular $\mathbf{z}_\varepsilon\not\in\mathcal
			C_\mathcal{T}(\eta)$ vertex for the $h$-method}
		\label{fig:Square_exact_h}
		
			\end{figure}
			Fig. \ref{fig:Square_exact_h} displays the total error
			$\|\nabla(\mathbf{u}-\mathbf{u}_\mathbf{S})\|_{\mathbb L^2(\Omega)} + \|p-p_M\|_{L^2(\Omega)}$ and the $L^2$ norm of the
			discrete divergence $\|\operatorname*{div}\mathbf{u}_\mathbf{S}\|_{L^2(\Omega)}$ when $\mathbf{z}_\varepsilon\not\in\mathcal{C}_\mathcal{T}(\eta)$ is treated as
			non-singular or as singular $\mathbf{z}_\varepsilon\in\mathcal{C}_\mathcal{T}(\eta)$ vertex.
			In the first case, our pressure-wired Stokes FEM is identical to the classical Scott-Vogelius FEM. 
			We observe that
			for low polynomial degrees $k$ and perturbations $\Theta(\mathbf{z}_{\varepsilon})>10^{-6}$, both variants
			converge optimally.
			However, for small perturbations $\Theta(\mathbf{z}_{\varepsilon})<10^{-6}$ (dotted), the solution to the classical Scott-Vogelius FEM is polluted by
			the high condition number of the algebraic system and only our modification that treats
			$\mathbf{z}_\varepsilon\in\mathcal{C}_\mathcal{T}(\eta)$ as \W{$\eta$-critical} converges at all and with optimal rate.
			The divergence of the classical Scott-Vogelius FEM vanishes up to rounding errors except for a similar pollution
			effect for $\Theta(\mathbf{z}_\varepsilon)<10^{-6}$. 
			We also observe that the $L^2$ norm of the divergence for the pressure-wired Stokes FEM is very small compared to the
			total error and also significantly smaller than the velocity error
			$\|\nabla(\mathbf{u}-\mathbf{u}_\mathbf{S})\|_{\mathbb L^2(\Omega)}$ (undisplayed) by a factor of
			$\Theta(\mathbf{z}_{\varepsilon})$.
			\newline\textbf{Exponential convergence of the $k$-method.}
			This benchmark monitors the behaviour of the $k$-method that successively increases the polynomial
			degree $k$. 
			The convergence history of the total error in Fig. \ref{fig:Square_exact_p} displays the same pollution effect for
			small perturbations $0<\varepsilon\ll1$ or higher polynomial degree $k$ when $\mathbf{z}_\varepsilon$ is
			treated as a non-singular vertex.
			All other graphs overlay and, in particular, the pressure-wired Stokes FEM converges exponentially.
			For the discrete divergence, we observe the same convergence up to a certain threshold when accumulated
			rounding errors dominate.
			We remark that a sophisticated choice of bases, e.g., from \cite{beuchler_new_2006}, could improve this
			threshold.
			However, this does not affect the high condition number caused by the singular mesh configuration that produces the
			pollution effect for
			the Scott-Vogelius FEM.

\bibliography{sn-bibliography-new}
\end{document}